\documentclass[VANCOUVER,STIX1COL]{WileyNJD-v2}

\usepackage{algorithm,algorithmicx,algpseudocode}
\usepackage[switch]{lineno}
\usepackage{subcaption}
\usepackage{diagbox}
\usepackage{xcolor}

\newlength{\commentindent}
\articletype{Article Type}%

\usepackage{cleveref}
\crefformat{subsection}{\S#2#1#3}
\crefformat{subsubsection}{\S#2#1#3}

\def\2nm#1{\|#1\|_2}

\def\1nm#1{\|#1\|_1}

\received{n/a}
\revised{n/a}
\accepted{n/a}
\begin{document}

\title{{Scaled ILU Smoothers for Navier-Stokes Pressure Projection}}

\author[1]{Stephen Thomas*}

\author[2]{Arielle K. Carr}

\author[3]{Paul Mullowney}

\author[4]{Katarzyna \'{S}wirydowicz}

\author[5]{Marcus S. Day}

\authormark{Thomas, Carr, Mullowney, \'{S}wirydowicz, Day}

\address[1]{\orgname{National Renewable Energy Laboratory}, \orgaddress{Golden, \state{CO}, \country{USA}}, \email{Stephen.Thomas@nrel.gov}}

\address[2]{\orgdiv{Department of Computer Science and Engineering}, \orgname{Lehigh University}, \orgaddress{Bethlehem, \state{PA}, \country{USA}}, \email{arg318@lehigh.edu}}

\address[3]{\orgname{AMD}, \orgaddress{Boulder, \state{CO}, \country{USA}}, \email{Paul.Mullowney@amd.com}}

\address[4]{\orgname{Pacific Northwest National Laboratory}, \orgaddress{Richland, \state{WA}, \country{USA}}, \email{kasia.swirydowicz@gmail.com}}

\address[5]{\orgname{National Renewable Energy Laboratory}, \orgaddress{Golden, \state{CO}, \country{USA}}, \email{Marcus.Day@nrel.gov}}

\corres{*Corresponding author}
 

\abstract[Summary]{
Incomplete LU (ILU) smoothers are effective in the algebraic multigrid 
(AMG) $V$-cycle for reducing high-frequency 
components of the error. However, the requisite direct triangular 
solves are comparatively slow on GPUs. Previous work has 
demonstrated the advantages of
Jacobi iteration as an alternative to direct solution of 
these systems.
Depending on the threshold and fill-level parameters chosen, 
the factors can be highly non-normal and Jacobi is unlikely 
to converge in a low number of iterations.  
We demonstrate that row scaling can reduce the departure 
from normality, allowing us to replace the inherently 
sequential solve with a rapidly converging Richardson iteration.
There are several advantages beyond the lower compute time. 
Scaling is performed locally for a diagonal block of the global matrix 
because it is applied directly to the factor. Further, an ILUT Schur 
complement smoother 
maintains a constant GMRES iteration count as the number of MPI ranks increases,
and thus parallel strong-scaling is improved. 
Our algorithms have been incorporated into hypre, 
and 
we demonstrate improved time to solution for {linear systems arising in the}
Nalu-Wind and PeleLM
pressure solvers. For large problem sizes, 
GMRES$+$AMG executes at least five 
times faster when using iterative triangular solves  compared with direct solves on massively-parallel GPUs.}

\keywords{iterative triangular solvers, GPU acceleration, ILU smoother, algebraic multigrid}

\jnlcitation{\cname{%
\author{Thomas S. J.},
\author{Carr, A. K.},
\author{Mullowney, P.} 
\author{\'{S}wirydowicz K.}, and
\author{Day, M.}} (\cyear{2022}), 
\ctitle{Scaled ILU Smoothers for Navier-Stokes Pressure Projection}, \cjournal{NAME}, \cvol{VOL}.}

\maketitle

\section{Introduction}\label{sec:intro}

{
We consider an incomplete LU (e.g., ILU$(k)$, ILUTP) \cite{Saad03} factorization
used within the Algebraic Multigrid (AMG) method for the low Mach Navier-Stokes
pressure solvers, PeleLM\cite{PeleLM} and Nalu-Wind\cite{Thomas2019}.  
Analogous geometric multigrid methods (GMG) are employed in the NekRS \cite{fischer2022nekrs,Nek5000Website}
and ExaDG \cite{arndt2020exadg} pressure solvers.
GMG operates on a hierarchy of nested grids with different resolutions and is 
effective for structured grid problems where the underlying mesh has a 
well-defined geometric structure. 
Our applications employ unstructured grids for which AMG 
is the appropriate choice.
When employing ILU as smoothers within AMG, 
fast and accurate solvers for the
resulting triangular systems are essential to achieve a low backward error.
However, directly solving the triangular systems in the AMG solve phase may
result in a significant performance bottleneck on massively parallel
architectures.  Iterative methods, like Richardson, offer an
efficient alternative for approximating the sparse triangular solution
\cite{Anzt2015, Chow2018}, but these methods may fail to converge in a
sufficiently small number of iterations when the $L$ or $U$ factor exhibit a
high degree of non-normality. This can occur when factorizing ill-conditioned
coefficient matrices\cite{ChowSaad1997} and is the case with the two applications considered.
When the triangular factors are close to normal, it
becomes possible to use the significantly faster sparse matrix-vector (SpMV)
products appearing in the Richardson iteration on GPU architectures. Thus,
effective iterative methods avoid a large departure from normality, dep($U$) or
dep($L$), of the triangular factors. 
The large departure from normality and high condition number $\kappa(B)$
are mitigated by applying row 
scaling to the triangular factors. Substantial acceleration is thus achieved in the AMG
solve phase for matrices of dimension larger than 10M exported from
PeleLM\cite{PeleLM} and Nalu-Wind\cite{Thomas2019} compared with direct solves
and iterative solves without scaling.}

{Linear systems of the form $Ax=b$, where $A$ is a sparse
symmetric $n\times n$ matrix, arise in ``projection'' methods for evolving
variable-density incompressible and reacting flows in the low Mach flow regime.
$A$ is highly ill-conditioned when using cut-cell approaches to complex
geometries, where non-covered cells that are cut by the domain boundary can
have arbitrarily small volumes and areas \cite{PeleLM}. For instance, the
matrices from the PeleLM combustion model in section \ref{sec:peleresults} have
singular values spanning sixteen orders of magnitude. Re-ordering
the unknowns into blocks is one technique for handling non-normal 
factors\cite{Chow2018}. However, a reordering algorithm such
as the reverse Cuthill-McKee algorithm (RCM) can increase the
computational cost.}  

{Equilibration (i.e., scaling) techniques are generally
designed to reduce the condition number of $A$ \cite{VanDerSluis1969,Bauer1963}
and can therefore improve the accuracy of the solution to the linear system,
particularly when using direct techniques\footnote{The instability of the
solution to a linear system with an ill-conditioned coefficient matrix is
well-known; see \cite[Chapter 4.2]{Stewart1998}}.  However, equilibration has
received considerably less attention when used in conjunction with AMG. 
Even when equilibration is not performed, the $L$ and $U$ factors can be
highly non-normal \cite{ChowSaad1997}. Indeed,
equilibration of ill-conditioned $A$ itself can also lead to highly non-normal
factors $L$ and $U$ when subsequently computing the ILU factorization of $A$.  
In either case, when employing ILU
smoothers that result in triangular factors with a large departure from
normality, the use of Richardson iterative methods is rendered ineffective for
approximately solving the triangular systems in the AMG solve phase. In these
cases where we cannot reasonably avoid a high departure from normality a
priori, and instead seek to reduce the departure from normality after the ILU
factorization but prior to the triangular solve. We propose a technique to
scale highly non-normal $U$ factors using row scaling, $A \approx LDU$, where
$D$ is a diagonal matrix extracted from $U$ by row scaling.\footnote{Note that
$D$ can also be factored out from $L$.} dep$(L)$ and $\kappa(L)$ remain modest, 
and thus we scale only the $U$ factor. However, our approach extends to 
scaling both factors if necessary. Scaling leads to a significant 
reduction in both the condition number
and departure from normality of the $U$ factor \cite{Henrici1962}.} 

{We also consider an ILUT Schur complement smoother and
demonstrate that this maintains a constant Krylov solver iteration count as the
number of MPI ranks increases. The Schur
complement system represents the interface degrees of freedom at subdomain
boundaries and are associated with the column indices corresponding to row
indices owned by other MPI ranks in the hypre ParCSR block partitioning of the
global matrix \cite{Falgout2004}.  A single GMRES iteration is employed to
solve the Schur complement system for the interface variables, followed by
back-substitution for the internal variables using iterative triangular solves.
A key observation is that the explicit residual computation $r^{(k)} = b -
Ax^{(k)}$ is not needed for this single GMRES iteration, resulting in
additional and significant computational cost saving.  A hierarchical basis
formulation of AMG based on the C-F block matrix partitioning and the Schur
complement has been considered previously\cite{Chow2004}, but can become
expensive as the number of nonzeros in the coefficient matrix increases.  A
similar increase in cost is observed in our experiments, along with an increase
in dep$(U)$ and $\|U_s\|_2$. In our approach, these
are mitigated by limiting the fill-in and with row scaling and
the Richardson iteration can again be effectively applied, resulting in the
comparatively inexpensive solution of the triangular systems for the internal
variables on the GPU.}
This paper is organized as follows.  Section 2 presents the low Mach
Navier-Stokes equations for the PeleLM and Nalu-Wind models.
{The AMG method} is reviewed in \cref{sec:AMG}, along with
smoothing techniques.  Departure from normality and row scaling are also
discussed.  {Equilibration of the triangular factors} for the
ILU smoother is presented in \cref{sec:ILUSmooth}.  Implementation of AMG
within the hypre \cite{Falgout2004} solver library is described in
\cref{sec:results}, where a performance model is provided, and results are
presented for linear systems from the PeleLM \cite{PeleLM} and Nalu-Wind
\cite{SC21} pressure solvers.  Conclusions are provided in section \cref{sec:concl},

{{\bf Contributions.} To our knowledge, scaling of the $L$ or
$U$ factors to reduce the condition number of the factors and thus mitigate
high non-normality is a novel approach when applied within an ILU smoother for
AMG.  Incorporating our scaling strategy, along with the ILUT Schur complement
smoother, into the AMG method applied to linear systems extracted from the
PeleLM nodal pressure projection solver\cite{PeleLM} and the Nalu-Wind pressure
continuity solver\cite{Thomas2019} facilitates a GMRES+AMG execution time of at
least five times faster when using iterative triangular solves on massively
parallel GPUs.  Scaling, in particular, reduces the condition number of $U$ and
dep$(U)$ significantly, facilitating the use of a truncated Neumann series
(Richardson) iteration to solve the triangular system. For
the applications considered, scaling $L$ was found to be
unnecessary, but our approach can be easily extended to cases when scaling is
needed for both triangular factors. Our
algorithms were incorporated into hypre\cite{Falgout2004}. We demonstrate the effectiveness and
efficiency of our approach for important problems arising in low Mach
Navier-Stokes pressure solvers solved using AMG.  }

\section{Low Mach Incompressible Navier-Stokes}

\subsection{PeleLM Combustion Model}

PeleLM is an adaptive mesh
low Mach number combustion code developed and supported under
DOE's Exascale Computing Project. PeleLM features
a variable-density projection scheme to ensure that the
velocity field used to advect the state satisfies an elliptic
divergence constraint. Physically, this constraint enforces a
consistently evolving flow with a spatially uniform
thermodynamic pressure across the domain. A key feature of
the model is that the discretization is based on a conservative 
embedded boundary approach to represent complex boundary shapes.
Intersections of domain boundaries with the underlying 
Cartesian grid can lead to arbitrarily small cell faces and volumes, 
which in turn can lead
to highly
ill-conditioned matrices representing the elliptic projection
operator. 

The low Mach number flow equations represent
the reacting Navier-Stokes flow equations in the low Mach number regime,
where the characteristic fluid velocity is small compared to the sound speed, 
and the effect of acoustic wave propagation is unimportant to the overall dynamics 
of the system. Accordingly, acoustic wave propagation can be mathematically removed 
from the equations of motion, allowing for a numerical time step based on an advective CFL condition. This leads to an increase in the allowable time step of order $1/M$ over an explicit, 
fully compressible method, where $M$ is the Mach number.  In this mathematical framework, the total 
pressure is decomposed into the sum of a spatially constant (ambient) thermodynamic pressure 
and a perturbational pressure, $\pi(x,t)$, that drives the flow.  
The set of conservation equations specialized to the low Mach number regime is a system of
partial differential equations with advection, diffusion, and reaction (ADR) processes that are constrained to evolve 
on the manifold of a spatially constant $P_0(t)$. Under suitable conditions, 
$\pi(x,t)/P_0(t) = \mathcal{O} (M^2)$. 

\[\begin{split}\frac{\partial (\rho \boldsymbol{u} )}{\partial t} +
\nabla \cdot \left(\rho  \boldsymbol{u} \otimes  \boldsymbol{u}  + \tau \right)
= -\nabla \pi + \rho \boldsymbol{F},\\
\frac{\partial (\rho Y_m)}{\partial t} +
\nabla \cdot \left( \rho Y_m \boldsymbol{u}
+ \boldsymbol{\mathcal{F}}_{m} \right)
= \rho \dot{\omega}_m,\\
\frac{ \partial (\rho h)}{ \partial t} +
\nabla \cdot \left( \rho h \boldsymbol{u}
+ \boldsymbol{\mathcal{Q}} \right) = 0 ,\end{split}\]
$\rho$ is the density, $\boldsymbol{u}$ is the velocity,
is the mass-weighted enthalpy,
$T$ is temperature, $Y_m$ is the mass fraction of species, and $\dot{\omega}_m$ 
is the molar production rate for species $m$. Additionally, $\tau$ is the stress tensor,  $\boldsymbol{\mathcal{Q}}$ 
is the heat flux and $\boldsymbol{\mathcal{F}}_m$ 
are the species diffusion fluxes. 
These transport fluxes require the evaluation of transport coefficients 
(e.g., the viscosity $\mu$, the conductivity, and the diffusivity matrix 
$D$, all of which are computed using the library EGLIB~\cite{EGLIB}). 
The momentum source,  $\boldsymbol{F}$,
is an external forcing term.
These evolution equations are supplemented by an equation of state 
for the thermodynamic pressure.  For example, the ideal gas law,
\[P_0(\rho,Y_m,T)=\frac{\rho \mathcal{R} T}{W}=\rho \mathcal{R} T
\sum_m \frac{Y_m}{W_m},\]
can be used, where $W_m$ and 
$W$ are the species and mean molecular weights, respectively.  
To close the system we also require a relationship between enthalpy, 
species, and temperature.  We adopt the definition used in the CHEMKIN
~\cite{CHEMKIN} standard, and have
$h=\sum_m Y_m h_m(T)$,
where $h_m$ is the species enthalpy.  Note that expressions for 
$h_m(T)$ incorporate the heat of formation for each species.
Neither species diffusion nor reaction redistribute the total mass, 
hence $\sum_m \boldsymbol{\mathcal{F}}_m = 0$ and 
$\sum_m \dot{\omega}_m = 0$. Thus, summing the species equations 
and using the definition 
$\sum_m Y_m = 1$, we obtain the continuity equation:
\[\frac{\partial \rho}{\partial t} + \nabla \cdot \rho \boldsymbol{u} = 0\]
This, together with the conservation equations form a 
differential-algebraic equation (DAE) system that describes 
an evolution equation subject to the constraint of spatially constant 
thermodynamic pressure.  A standard approach to 
attacking such a system computationally is to differentiate the 
constraint until it can be recast as an initial value problem.  
Following this procedure, set the thermodynamic pressure 
constant in the frame of the fluid,
$DP_0/Dt = 0$
and observe that if the initial conditions satisfy the constraint, 
this expression will enforce that it continues to 
satisfy the  constraint over all time.  
Expanding this expression via the chain rule and continuity:
\[\nabla \cdot \boldsymbol{u} = \frac{1}{T}\frac{DT}{Dt}
+ W \sum_m \frac{1}{W_m} \frac{DY_m}{Dt} = S.\]
Thus the constraint here takes the form of a condition on the divergence of the flow. 
Note that the actual expressions will depend upon the chosen 
models for evaluating the transport fluxes.

There are three different types of linear solves required to advance the 
velocity field. The first is the marker-and-cell (MAC) solve in order to obtain face-centered
velocities used to compute advective fluxes. The second is the multi-component
cell-centered solver that is used to obtain the provisional new-time velocities.
Finally, a nodal solver is used to project these in order that 
they satisfy the constraint.
We project the new-time velocity by solving the elliptic equation,
\begin{equation}
L\:\phi = D \: \left[ \: 
\boldsymbol{u}^{n+1,*} + \frac{\delta t}{\rho^{n+1/2}}\:G\:\pi^{n-1/2}
\: \right] - \hat{S}^{n+1}
\label{eq:nodalP}
\end{equation}
for nodal values of $\phi$ and time index $n$. Here, $L$ represents a Laplacian of nodal 
data, computed using the standard bilinear finite-element approximation to 
$\nabla\:(1/\rho^{n+1/2})\:\nabla$. Also $D$ is a discrete second-order 
operator that approximates the divergence at nodes from cell-centered data 
and $G$ approximates a cell-centered gradient from nodal data. 
Nodal values for $\hat{S}^{n+1}$  required for this equation are obtained 
by interpolating the cell-centered values. Finally, we determine the
new-time cell-centered velocity field using
\begin{equation}
 \boldsymbol{u}^{n+1} =  \boldsymbol{u}^{n+1,*} -
\frac{\delta t}{\rho^{n+1/2}}\:G\:(\:\phi - \pi^{n-1/2}\:)
\end{equation}

\subsection{ExaWind Wind-Turbine Model}

The ExaWind ECP project aims to simulate the atmospheric boundary layer
air flow through an entire wind farm on next-generation exascale-class computers.
The primary physics codes in the ExaWind simulation environment are
Nalu-Wind and AMR-Wind. Nalu-Wind and AMR-Wind are finite-volume-based CFD 
codes for the incompressible-flow Navier-Stokes governing equations. Nalu-Wind is an
unstructured-grid solver that resolves the complex geometry of wind turbine 
blades and thin blade boundary layers.  AMR-Wind is a block-structured-grid 
solver with adaptive mesh refinement (AMR) capabilities that captures the 
background turbulent atmospheric flow and turbine wakes.  Nalu-Wind and 
AMR-Wind models are coupled through overset meshes. The equations consist 
of the mass-continuity equation for pressure and Helmholtz-type
equations for transport of momentum and other scalars (e.g. those for turbulence
models). For Nalu-Wind, simulation times are dominated by linear-system setup
and solution of the continuity and momentum equations. Both PeleLM and 
AMR-Wind are built on the AMReX software stack \cite{Zhang2019}, with
hierarchical block structured meshes, and employ
geometric multigrid as the primary solver,
however, they both have the option of using the hypre library \cite{Falgout2002} 
to access alternative solvers, such as algebraic multigrid.

An implicit BDF time integrator was employed with an adaptive (variable) time-step
that can increase towards a target Courant-Friedrichs-Lewy (CFL) number.  The momentum 
and pressure solutions obtained during time stepping are calculated by an incremental
approximate pressure projection algorithm. The momentum and pressure equations
are segregated and solved sequentially with implicit advection/diffusion.
To describe the approximate pressure projection algorithm, consider the block
(indefinite, saddle-point) matrix form of the discrete equations
\begin{equation}
\left[
\begin{array}{cc}
F  &  G   \\
D  &  0
\end{array}
\right]
\left[
\begin{array}{c}
u   \\
p
\end{array}
\right] =
\left[
\begin{array}{c}
f   \\
0
\end{array}
\right] \,,
\end{equation}
where the matrix $F$ will depend on the predicted, current and earlier time-levels. $F$ contains
discrete contributions to the momentum equations from the time derivative, diffusion,
and linearized advection terms. $F = I/\Delta t + \mu L + N$,
where $L$ is the Laplacian and $N$ linearized advection.
The discrete gradient and divergence matrices are $G$ and $D$ respectively.
The vector $f$ contains the additional terms for the
momentum equations, e.g., body force terms, lagged stress tensor terms, etc.
The right-hand side contains the appropriate terms, e.g., for a 
non-solenoidal velocity field. For more details, see Thomas et al.~\cite{Thomas2019}.
In order to derive a projection scheme, consider the time-split system of equations,
\begin{equation}
\left[
\begin{array}{cc}
F  &  G   \\
D  &  0
\end{array}
\right]
\left[
\begin{array}{c}
u^{n+1}   \\
\Delta p^{n+1}
\end{array}
\right]
+
\left[
\begin{array}{cc}
I  &  G   \\
D  &  0
\end{array}
\right]
\left[
\begin{array}{c}
0         \\
p^n
\end{array}
\right] =
\left[
\begin{array}{c}
f   \\
0
\end{array}
\right] \,,
\end{equation}
where $\Delta p^{n+1} = p^{n+1} - p^n$. Consider the block factorization of
the matrix
\begin{equation}
\left[
\begin{array}{cc}
F  &  G   \\
D  &  0
\end{array}
\right]
=
\left[
\begin{array}{cc}
F  &  0   \\
D  &  -DF^{-1}G
\end{array}
\right]
\left[
\begin{array}{cc}
I  &  F^{-1}G   \\
0  &  I
\end{array}
\right]  \,.
\end{equation}
Inversion of $F$ to form the Schur complement matrix $M=-DF^{-1}G$ would
be costly. The splittings
approximate the inverse with the diagonal matrix $[ {\rm diag}(F)]^{-1}$.
Nalu-Wind employs an approximate projection scheme that introduces an auxiliary
projection time-scale, determined by the
factor $B_2 = (\tau/\rho)\:I$.
The factor $B_1=-\tau L$ defines the linear system for pressure and
approximates $M$.  $L$ is the Laplacian matrix obtained from the discrete
form of the Gauss divergence theorem
%
%

These matrices are introduced into an approximate block factorization as follows
\[
\left[
\begin{array}{cc}
F  &  0   \\
D  &  B_1
\end{array}
\right]
\left[
\begin{array}{cc}
I  &  B_2G   \\
0  &  I
\end{array}
\right]
\left[
\begin{array}{c}
u^{n+1}   \\
\Delta p^{n+1}
\end{array}
\right] +
\left[
\begin{array}{cc}
I      &  0   \\
DB_2 &  -B_1
\end{array}
\right]
\left[
\begin{array}{cc}
I  &  G   \\
0  &  I
\end{array}
\right]
\left[
\begin{array}{c}
0         \\
p^n
\end{array}
\right] =
\left[
\begin{array}{c}
f   \\
0
\end{array}
\right] \,.
\]
%
The time scale $\tau = \Delta t$ is chosen for stabilization.
The equations are solved for $\Delta \hat{u} = \hat{u} - u^n$, and $\Delta p^{n+1}$
at each outer nonlinear iteration of the time step
\begin{eqnarray}
F\Delta \hat{u}         & = & f - Gp^n - Fu^n \,,
\label{eq:momen} \\
\tau\:L\Delta p^{n+1}   & = & D\: ( \:\rho\hat{u} + \tau\:G\:p^n\:) - \tau\: L\:p^n   \,,
\label{eq:press} \\
u^{n+1}                 & = & \hat{u} - \frac{\tau}{\rho} \: G \Delta p^{n+1} \,.
\label{eq:correct}
\end{eqnarray}
The momentum (\ref{eq:momen}), and pressure-continuity (\ref{eq:press})
equations are solved separately within a fixed-point iteration for the
incompressible system of equations. An iteration predicts a
velocity field that is not necessarily divergence free and then projects that field to a divergence
free sub-space.  The matrix $F$ and the right-hand side of the momentum
equation are functions of the solution at time step $n+1$ due to the choice of
an implicit BDF time integrator. At each iteration, a better estimate for the
solution is computed at time step $n+1$ and hence a better estimate for $F$.

The pressure projection equations (\ref{eq:nodalP}) and (\ref{eq:press})
are  solved in discrete form with the GMRES+AMG iterative method
described in the following sections. AMG is an optimal
solver because the amount of computational work per degree
of freedom in the resulting linear system remains constant as the
problem size increases.


\section{Algebraic Multigrid}\label{sec:AMG}
The AMG implementation to solve the linear system $Ax=b$ consists of two steps: the 
setup and solve phases. In the setup phase, the method constructs {\it prolongation} 
and {\it restriction} operators to transfer field data between coarse and fine 
grids.  The application of these transfer operators leads to a sequence, or 
{\it hierarchy}, of successively lower dimension matrices, denoted 
$A_{l} \in \mathbb{C}^{m_l\times m_l}$, 
$l=0,1,\dots m$ where $A_l=R_l \: A_{l-1}\: P_l$, $m_l < m_{l-1}$, and $m_0 = n$. 
In the Galerkin formulation of 
AMG, $P_l$ is a rectangular matrix with dimensions $m_{l-1} \times
m_l$ also referred to as the {\it interpolant} and $R_l = P_l^T$ . 
Once the transfer operators are determined, the
coarse-matrix representations are computed through
sparse triple-product matrix-matrix multiplication. 

In classical Ruge-St\"{u}ben AMG \cite{Ruge1987}, 
the {\it strength of connection threshold}
determines whether points in the fine mesh are retained when creating a 
set of coarse points: The point $j$ {\it strongly influences} the point 
$i$ if and only if
\begin{equation*}
|a_{ij}| \ge \theta \max_{k\ne i}\:|\:a_{ik}\:| \,,
\end{equation*}
where $\theta$ is the strength of connection threshold, 
$0 < \theta \le 1$. The selected 
coarse points are retained at the next coarser level, and the remaining 
fine points are dropped. More information on the strong connection 
threshold to determine the set of coarse grid points as well as
how to form the transfer operators are found in \cite{Ruge1987}. 
Let $C_l$ and $F_l$ be the coarse and fine points selected at
level $l$, and let $m_l$ be the number of grid points at level $l$.
Then, $m_l = |C_l| + |F_l|$ and $m_{l+1} = |C_l|$. Here, the coarsening is
performed row-wise by interpolating between coarse and fine levels and generally 
attempts to fulfill two contradictory criteria.
In order to ensure that a chosen interpolation scheme is well-defined and of
good quality, some close neighborhood of each fine point must contain a
sufficient amount of coarse points to interpolate from. Hence, the set of coarse
points must be rich enough {but should also} be
sufficiently small in order to achieve a reasonable coarsening rate.
Because the size of the linear systems decreases 
on each coarser level, the interpolation should lead to a sufficient 
reduction in the number of non-zeros at each level of the hierarchy 
compared with the number of non-zeros in the coefficient matrix for the original 
linear system.\footnote{While the coarser matrices are 
technically more dense - as a ratio of non-zeros to the size of the matrix - 
these matrices have fewer rows and columns.} 

In the solve phase, the method employs the {\it V-cycle}, which 
{is comprised of} a relaxation (or smoothing) iteration 
coupled with a coarse grid correction.\footnote{ The V-cycle is the simplest 
complete AMG cycle.  Other processes may be used in place of the V-cycle, 
such as the W- or F-cycles \cite{Ruge1987}.} Beginning at the finest level, 
the method moves to the next coarser level by first performing a small number 
of {\it pre-smoothing} iterations to 
the solution of $A_l x_l = b_l$ at the $l^{th}$ level of the $V$-cycle.  
It then computes the residual $r_l = b_l - A_l x_l$ and applies the 
restriction operator as $b_{l+1} = R_l r_l$ to move to the next coarser
level. This process is repeated until the coarsest level is reached.
The equations at the coarsest level are either solved
using a direct solver, or using an iterative scheme 
if $A_m$ is singular.  {Here, the method interpolates the coarse
solution as a correction to the solution on the next finer level using prolongation, 
$x_l = x_l + P_k\:x_{l+1}$}, 
followed by a small number 
of {\it post-smoothing} iterations.  
AMG methods are optimal for certain linear systems,
(i.e., constant work per degree of freedom in $A_0$) through complementary
error reductions by the smoother and solution corrections propagated from
coarser levels. Consult Algorithm \ref{fig:mgcycle} for a description 
of the multigrid $V$-cycle. 

\begin{algorithm}
\centering
\begin{algorithmic}[0]
  \State{//Solve $Ax=b$}
  \State{Set $x=0$}
  \State{Set $\nu=1$ for $V$-cycle}
  \State{call Multilevel($A, b, x, 0, \nu$)}
  \State{}
  \Function{Multilevel}{$A_l$, $b$, $x$, $l$, $\nu$}
    \State{// Solve $A_l x = b$ ($l$ is current grid level)}
    \State{// Pre smoothing step }
    \State $ x = S^{1}_l (A_l, b, x)$
      \If{$(k \ne {m})$}
        \State{$r_{l+1} = R_l (b - A_l x )$}
        \State{$A_{l+1} = R_l A_l P_l$}
        \State{$v = 0$}
        \For{$i = 1\dots\nu$}
          \State{}\Call{Multilevel}{$A_{l+1}$, $r_{l+1}$, $v$, $l+1$, $\nu$}
        \EndFor
        \State{$ x = x + P_{l} v$}
        \State{// Post smoothing step }
        \State{$ x = S^{2}_l (A_l, b, x )$}
      \EndIf
  \EndFunction
\end{algorithmic}
\caption{\label{fig:mgcycle} Multigrid single-cycle ($V$-cycle) algorithm 
for solving $Ax=b$.
}
\end{algorithm}

\subsection{Smoothing}\label{sec:ErrorSmooth}
Smoothers, {or relaxation schemes,} are generally implemented 
as an inexpensive iterative method such as Gauss-Seidel, 
Jacobi, or incomplete factorization, {and} rapidly reduce 
high-frequency components of the error by approximately solving the system of 
equations. {Given an approximate solution, $x_l$, to the system 
$A_lx_l = b_l$ at the $l^{th}$ level of the $V$-cycle, the general form of a 
residual smoothing technique can be expressed as $x^{(k+1)}_k = x^{(k)}_l + Sr^{(k)}_l$, 
where $x^{(k+1)}_l$ is the updated solution after smoothing, $S$ is the 
smoothing or damping factor, and $r^{(k)}_l = b_l-Ax^{(k)}_l$ is the residual.} 
When the remaining low-frequency error is restricted it then becomes higher frequency 
on the {next,} coarser level. A smoothing technique is employed
{until reaching the coarsest level of the $V$-cycle. Dropping the 
subscript $l$ for ease of notation,} the general form of a relaxation scheme 
for $Ax=b$ is given by
\begin{equation}\label{eq:splitting}
Mx^{(k+1)} = Nx^{(k)} + b
\end{equation}
where $A = M - N$ is a {\it matrix splitting}, {and $x^{(l+1)}$ 
represents the approximate solution at step $l+1$ of the iterative method, and 
$r^{(l)}$ is the residual at the $l^{th}$ step.} {For instance, the} 
Gauss-Seidel iteration
is based upon the splitting $M = D + L$, and $N = -U$, where $L$ is the strictly
lower and $U$ the upper triangular {part of the matrix $A$}. 
The inverse of the matrix $M$ is not formed, but rather direct triangular solvers
are typically employed. However, as {is the focus of this paper}, 
these are relatively slow on GPU architectures. {Adopting the notation 
used in \cite{Anzt2015}, the Jacobi iteration, is often written in the compact form 
\begin{equation}\label{eq:compact}
x^{(k+1)} = G\:x^{(k)} + D^{-1}b, 
\end{equation}
with the regular splitting $A = M - N$,
$M = D$ and $N = D - A$, and $G = I-D^{-1}A$.\footnote{{Note that 
here we are referring to Jacobi in the context of AMG smoothers. Later, we will refer 
to it as an iterative solver for the $L$ and $U$ systems that result from employing 
the ILU smoother.}}} 

Polynomial smoothers rely on fast {computation of the SpMV products}. 
A polynomial type smoother \cite{SC21} is derived from the
iterative solution of the triangular system for $(D + L)$ in Gauss-Seidel 
relaxation and then used to solve the linear system, $Ax = b$, where $D$ is 
{again} the diagonal of $A$. An
alternate formulation is to replace $(D + L)^{-1}$ with 
$(I+D^{-1}L)^{-1} D^{-1}$ in the preconditioned iteration, and
replace the matrix inverse with a truncated Neumann series
\[
x^{(k+1)} = x^{(k)} +  \sum_{j=0}^{p}
(-D^{-1} L)^j D^{-1}\:r^{(k)}, 
\]
for $p < n$. 


An ILU smoother is well-suited to handle highly varying matrix coefficients
or anisotropic problems and is a generalization of
Gauss-Seidel, where the diagonal matrix $D$ represents row scaling of one 
(or both) of the triangular factors {and $M = LDU$.  
Then we can write (\ref{eq:splitting}) as} 
\begin{equation}\label{eq:ILUsplit}
    x^{(k+1)} = x^{(k)} + M^{-1}\: (\: b - A\:x^{(k)} \:)
\end{equation}
An ILU factorization {is particularly useful for sparse matrices, 
maintaining the sparsity pattern of the original matrix.} It can be split 
into symbolic and numeric phases and if the
sparsity pattern of the matrix $A$ does not change on the finest level,
then the symbolic phase of the factorization is often reused. For example,
the sparsity pattern of the pressure matrix in the PeleLM model does not
change during the fluid integration step, but could change between time steps
as a result of a mesh refinement or regridding operation.
This permits re-use of factorizations in the numeric phase 
in order to save computational time and avoid releasing and re-allocating storage for
the $L$ and $U$ factors. This optimization opportunity exists only for the finest level 
in the multigrid hierarchy because
the data-dependent coarsening algorithm may change the sparsity pattern on other levels.

An ILUT Schur complement smoother is also considered in\cite{Saad03} as part of a 
hierarchical basis formulation of AMG.  Following Saad \cite{Saad03}, 
Xu ~\cite{XuLiKuffuor2020} and Falgout et al.~\cite{Falgout2021},
a Schur complement preconditioner {can be derived} for the 
partitioned linear system
\[
A\left[
\begin{array}{c}
x \\
y
\end{array}
\right]
= 
\left[
\begin{array}{c}
f \\
g
\end{array}
\right].
\]
Consider the block $A = LU$ factorization of the coefficient matrix $A$
\[
A = \left[
\begin{array}{cc}
B & E \\
F & C
\end{array}
\right] =
\left[
\begin{array}{cc}
I & 0 \\
FB^{-1} & I
\end{array}
\right]
\left[
\begin{array}{cc}
B & E \\
0 & S
\end{array}
\right].
\]
The block matrix $B$ is associated with the diagonal block (sub-domain) of
the global matrix distributed across MPI ranks.
The Schur complement is $S = C - FB^{-1}E$, and the
reduced system for the interface variables, $y$, is given by
\begin{equation}
S\:y = g - F\:B^{-1}\:f
\label{eq:schur}
\end{equation}
Then the internal, or local, variables represented by $x$ are obtained by back-substitution
according to the expression
\[
x = B^{-1}\: (\: f - E\:y\:)
\]
An ILUT Schur complement smoother for one level of the $V$-cycle in
hypre is implemented as a
single iteration of a GMRES solver for the global interface system (\ref{eq:schur}). 
The local systems involving $B^{-1}$ are then solved by computing an ILUT
factorization of the matrix $B\approx LDU$. 

{We stress that the smoother employed on each level of the $V$-cycle 
does not have to be of the same type. Considering the aforementioned optimization of 
the ILU smoother on the finest level, this provides motivation for a mixed $V$-cycle 
with ILU smoothing only on the finest level for computational efficiency. However, 
even in the case of a mixed technique, when employing ILU smoothers on any level or 
all of them, care must be taken to ensure efficient solution of the resulting triangular 
systems. In Section \ref{sec:itertrisolve} we provide background on bottlenecks imposed 
by ill-conditioning and non-normality in ILU factorizations when performing iterative 
triangular solves, and in Section \ref{sec:ILUSmooth} we propose our equilibration 
technique for mitigating these problematic characteristics, facilitating the use of 
fast iterative methods for solving the resulting triangular systems on GPUs.}


\section{Iterative Triangular Solves}\label{sec:itertrisolve}

{A sparse triangular solver is a critical kernel in many 
scientific computing simulations, and significant efforts have 
been devoted to improving 
the performance of a general-purpose sparse triangular solver for 
GPUs~\cite{Li:2020}. For example, the traditional parallel algorithm 
is based upon level-set scheduling, derived from the sparsity structure 
of the triangular matrix.  Independent computations 
proceed within each level of the elimination tree.
Overall, the sparsity pattern of the matrix can result in an
extremely deep and narrow tree, thereby limiting the amount of available
parallelism for the solver to fully utilize many-core architectures,
especially when compared to an SpMV sparse matrix-vector product.
Multi-color reordering has been proposed as an alternative, however, this
approach can adversely impact the convergence rate of the solver.
In other words, there is little (or no) parallel work 
at most levels of the tree, and thus level scheduling does not provide
significant speed-up, if any at all.
Adopting an iterative approach leverages the speed of sparse 
matrix-vector products on GPUs. However, to facilitate fast convergence of the 
iterative triangular solves, we must address the concern of non-normality naturally 
introduced by the ILU factorization and exacerbated by poor choices in parameters.}

\subsection{{Convergence and Non-Normality}}\label{sec:nonnormal}
A {\it normal matrix} $A \in \mathbb{C}^{n\times n}$ satisfies 
$A^*A = AA^*$, and this property is referred to 
as {\it normality} throughout this paper.  
{Intuitively}, a {\it non-normal} matrix {can be} defined 
in terms of the difference between $A^*A$ and $AA^*$.  
In the current paper,  Henrici's definition of the 
{\it departure from normality} of a matrix  is employed
\begin{equation}\label{eq:depB}
    \text{dep}(A) = \sqrt{\|A\|_F^2-\|\Lambda\|_F^2},
\end{equation}
where $\Lambda\in\mathbb{C}^{n \times n}$ is the diagonal matrix containing the 
eigenvalues of $A$ \cite{Henrici1962}. 
Further 
information on metrics and bounds describing normality of matrices is found in
\cite{Ipsen1998,Henrici1962,Elsner1987} and references therein.  

In general, $\text{dep}(L)$ remains modest for the applications we consider in 
this paper, and the ILU factorization computes $L$ such that $L = I + L_s$. Thus, 
scaling is not applied to $L$ for these applications.  Again, we note that our 
methods easily extend to those applications that require scaling of both triangular 
factors. However, the same observations are not true for {dep}$(U)$. 
In particular, scaling is necessary to produce $U = I + U_s$, with unit diagonal, 
where $U_s$ is strictly upper triangular. Going forward, we will focus our discussion 
and notation on $U$.

{It follows directly from the definition of normality that an upper 
(or lower) triangular matrix cannot be normal unless it is a diagonal matrix 
(see \cite[Lemma 1.13]{Saad03} for a proof).  Therefore, some departure
from normality is expected in the $L$ and $U$.  However, if the 
departure from normality is too great, the iterations may diverge 
\cite{Anzt2015,Eiermann1993}.} When the number of non-zeros in the factors is
limited by larger drop tolerance and smaller fill-in levels 
{in the incomplete factorization process}, the number of
non-zeros in $U_s$ will decrease, and thus $\|U_s\|_F$ will become smaller.  
At the extreme, when fill-in is not allowed
beyond the diagonal, then $U_s = 0$, and 
\begin{equation}\label{eq:diagDep}
  \text{dep}(U) = \sqrt{\|I+U_s\|_F^2 - n} = \sqrt{\|I\|_F^2 - n} = 0.   
\end{equation} 
Thus, dramatically restricting the number of non-zeros in the factors is one 
approach for bounding $\text{dep}(U)$ a priori. However, imposing 
{constraints that result in factors that are too sparse} generally 
produces an ILU factorization that is {also} too inaccurate 
(e.g., $\|A-LU\|_F = \gamma \gg 0$) to be useful. In other words, if $LU$ is a poor 
approximation to $A$, we cannot reasonably expect an iterative solver to converge in 
few enough iterations to justify its use as a smoother. More generally, tuning the 
ILU parameters results in conservative or small $\text{nnz}(U)$ (and thus conservative 
$\text{nnz}(U_s)$), $\text{dep}(U)$ cannot grow too large.

\subsection{{Effects of ILU Parameter Choice on (Non-)Normality}}\label{sec:paramAnaly}
While the existing literature on the choice of parameters for 
ILU smoothers is limited, the work in \cite{ChowSaad1997} provides 
a framework for studying the effect of parameter choices.  
Here, the relationship between these choices 
and the resulting sizes of $\text{dep}(U)$ is briefly examined, 
and subsequently the number of Jacobi iterations required. 
For a more in depth analysis of parameter choices, we refer the reader to 
\cite{ChowSaad1997,Thomas2021Neumann}.
The potential inaccuracy when the $L$ and $U$ have high condition numbers 
is first examined. When employing a threshold parameter (or drop tolerance) 
to limit the number of non-zeros in a row (or column) of the 
factors, the factorization of a symmetric matrix could 
be highly nonsymmetric \cite{ChowSaad1997}.  
One observable indicator is the vertical striping in the sparsity pattern 
of $L+U$, which signifies orders of magnitude 
difference in the entries of a row (or column) of the 
coefficient matrix $A$ and is associated with ill-conditioning 
\cite{ChowSaad1997}.

\begin{figure}[hh]
\captionsetup{font=normalsize}
\captionsetup[subfigure]{justification=centering}
\begin{center} 
\begin{subfigure}{.23\textwidth}
\captionsetup{font=normalsize}
\includegraphics[width=\linewidth]{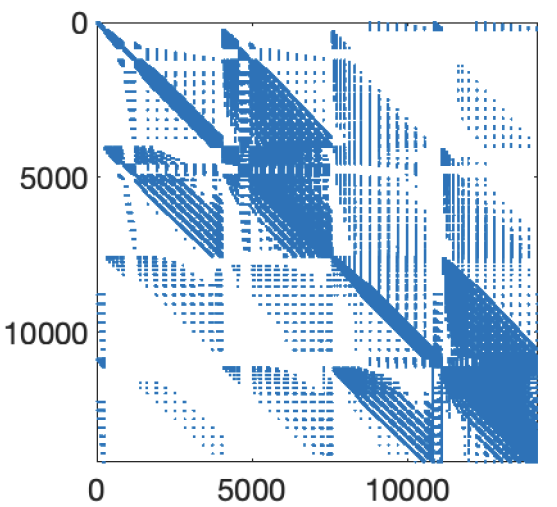}
\subcaption{First Level. \\$\kappa_2(U) =$ 1.29e13, \\$\kappa_2(L) =$ 1.95e8.}
\label{fig:LU1_N14k}
\end{subfigure}
\begin{subfigure}{.22\textwidth}
\captionsetup{font=normalsize}
\includegraphics[width=\linewidth]{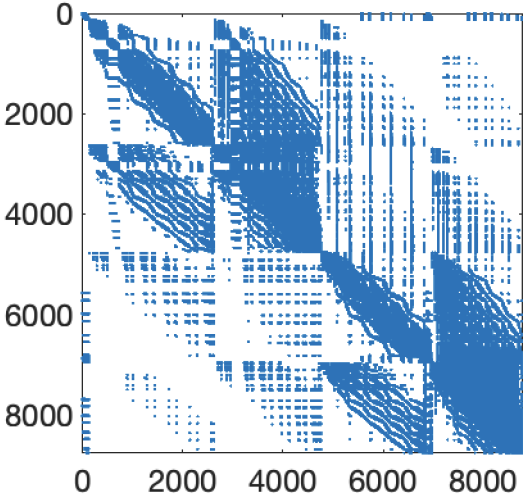}
\subcaption{Second Level.\\ $\kappa_2(U) =$ 1.25e13,\\ $\kappa_2(L) =$ 1.27e8.}
\label{fig:LU2_N14k}
\end{subfigure}
\begin{subfigure}{.22\textwidth}
\captionsetup{font=normalsize}
\includegraphics[width=\linewidth]{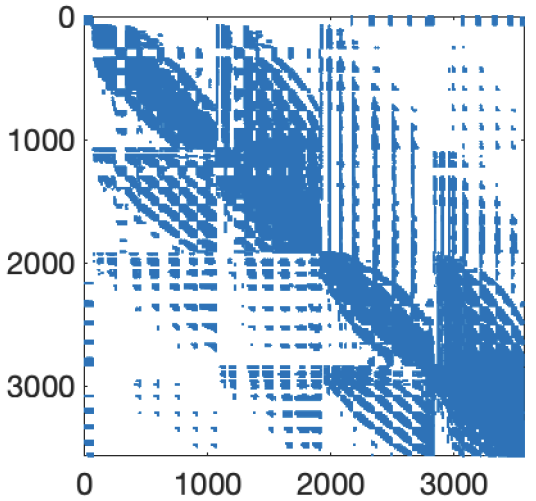}
\subcaption{Third Level.\\ $\kappa_2(U) =$ 1.01e14,\\ $\kappa_2(L) =$ 1.20e8.}
\label{fig:LU3_N14k}
\end{subfigure}
\begin{subfigure}{.22\textwidth}
\captionsetup{font=normalsize}
\includegraphics[width=\linewidth]{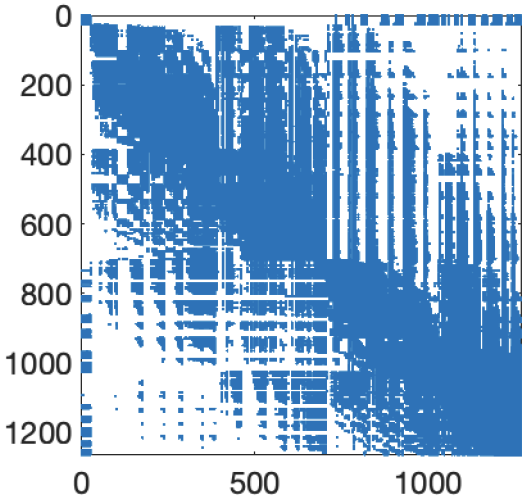}
\subcaption{Fourth Level.\\ $\kappa_2(U) =$ 6.53e16,\\ $\kappa_2(L) =$ 2.80e7.}
\label{fig:LU4_N14k}
\end{subfigure}
\end{center}
\caption{Non-zero patterns of $L+U$ for matrix size $N=14186$ for the first 
four levels using AMGToolbox.  Drop tolerance is set to 1.e$-15$ 
and fill limit per row set to 200.}
\label{fig:LU14k}
\end{figure}

\begin{figure}[hh]
\captionsetup{font=normalsize}
\captionsetup[subfigure]{justification=centering}
\begin{center} 
\begin{subfigure}{.23\textwidth}
\captionsetup{font=normalsize}
\includegraphics[width=\linewidth]{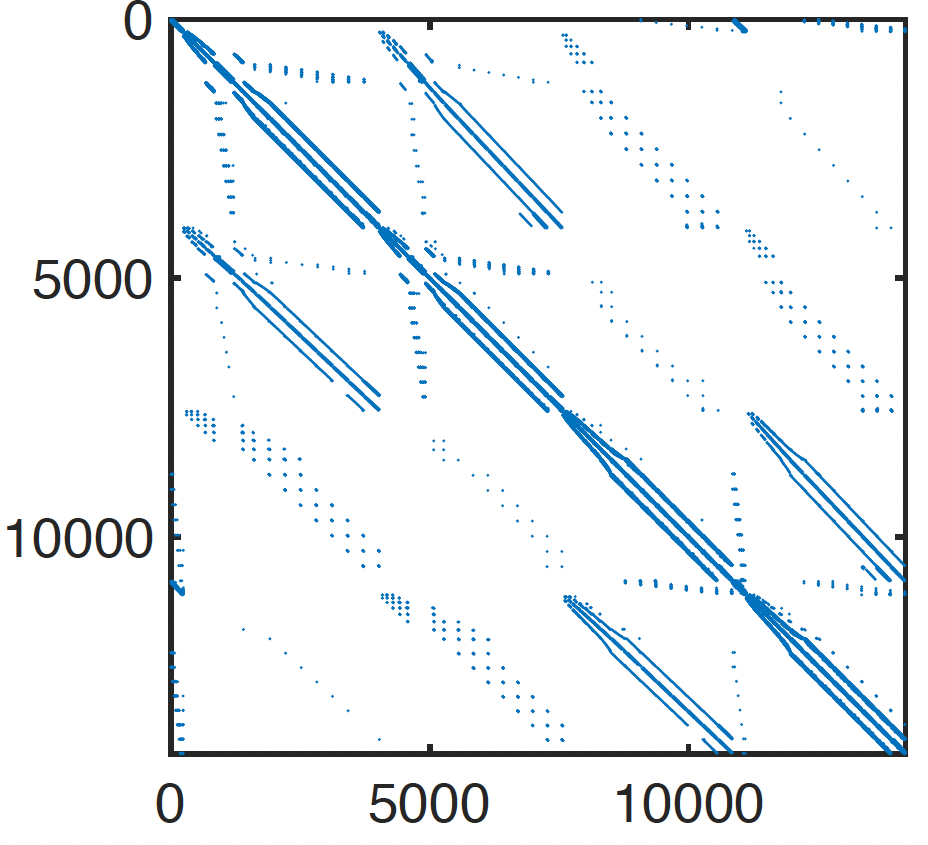}
\subcaption{First Level.\\ $\kappa_2(U) =$ 6.81e12,\\ $\kappa_2(L) =$ 1.94e8.}
\label{fig:LU1_N14k_fill}
\end{subfigure}
\begin{subfigure}{.22\textwidth}
\captionsetup{font=normalsize}
\includegraphics[width=\linewidth]{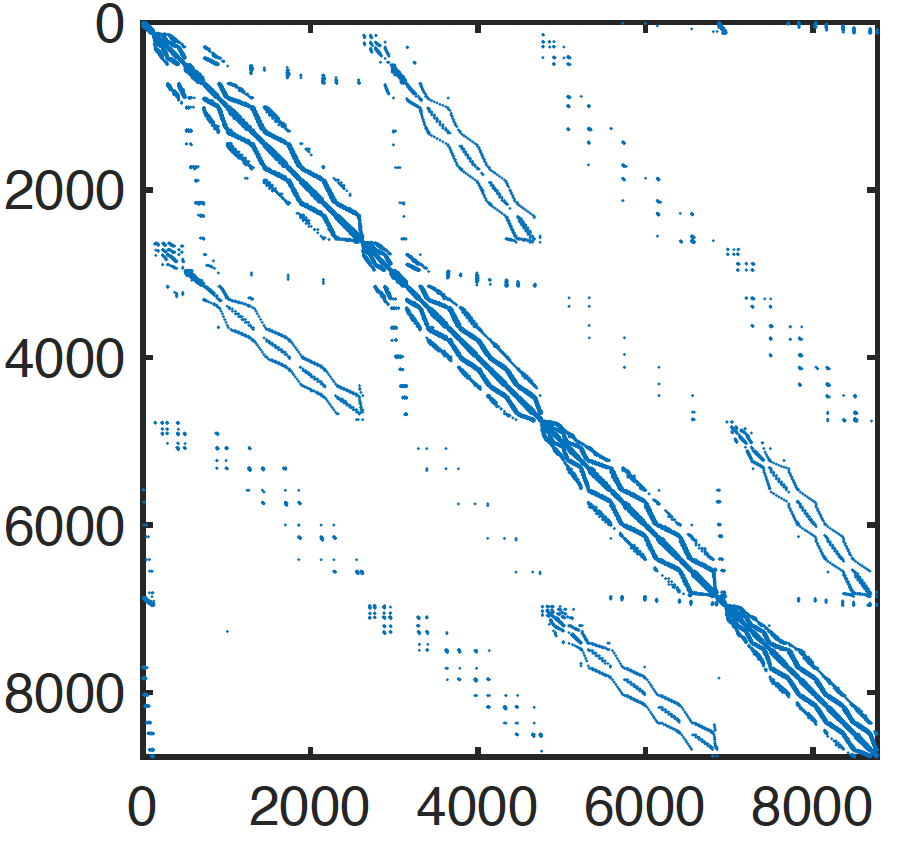}
\subcaption{Second Level.\\ $\kappa_2(U) =$ 1.00e13,\\ $\kappa_2(L) =$ 1.24e8.}
\label{fig:LU2_N14k_fill}
\end{subfigure}
\begin{subfigure}{.22\textwidth}
\captionsetup{font=normalsize}
\includegraphics[width=\linewidth]{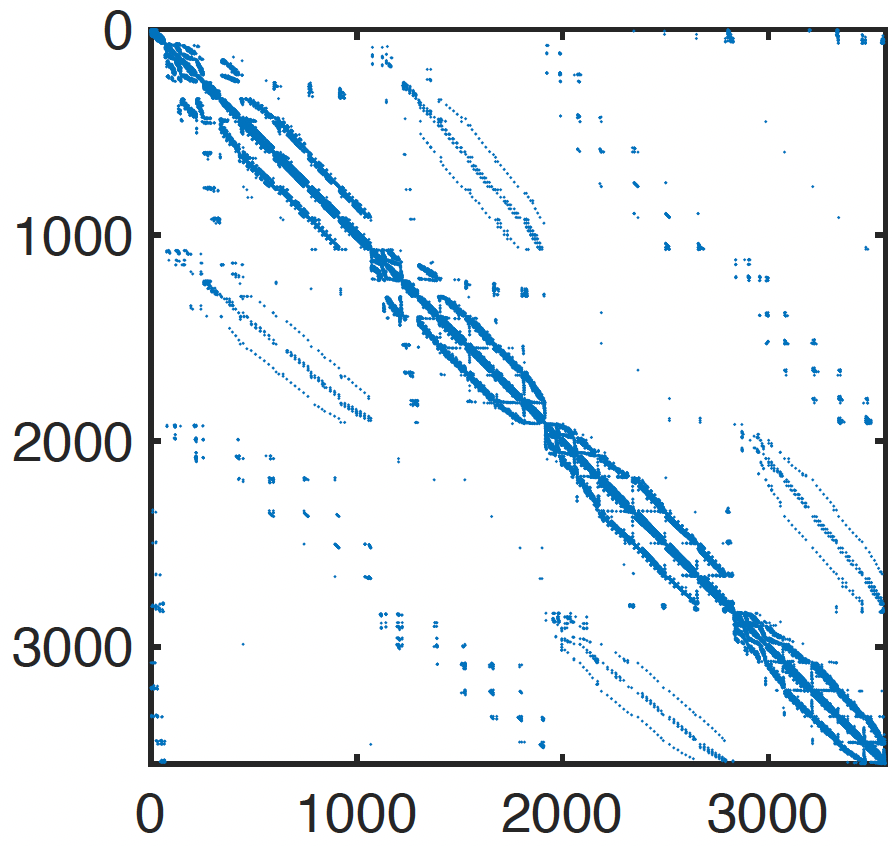}
\subcaption{Third Level.\\ $\kappa_2(U) =$ 9.97e13,\\ $\kappa_2(L) =$ 1.22e8.}
\label{fig:LU3_N14k_fill}
\end{subfigure}
\begin{subfigure}{.22\textwidth}
\captionsetup{font=normalsize}
\includegraphics[width=\linewidth]{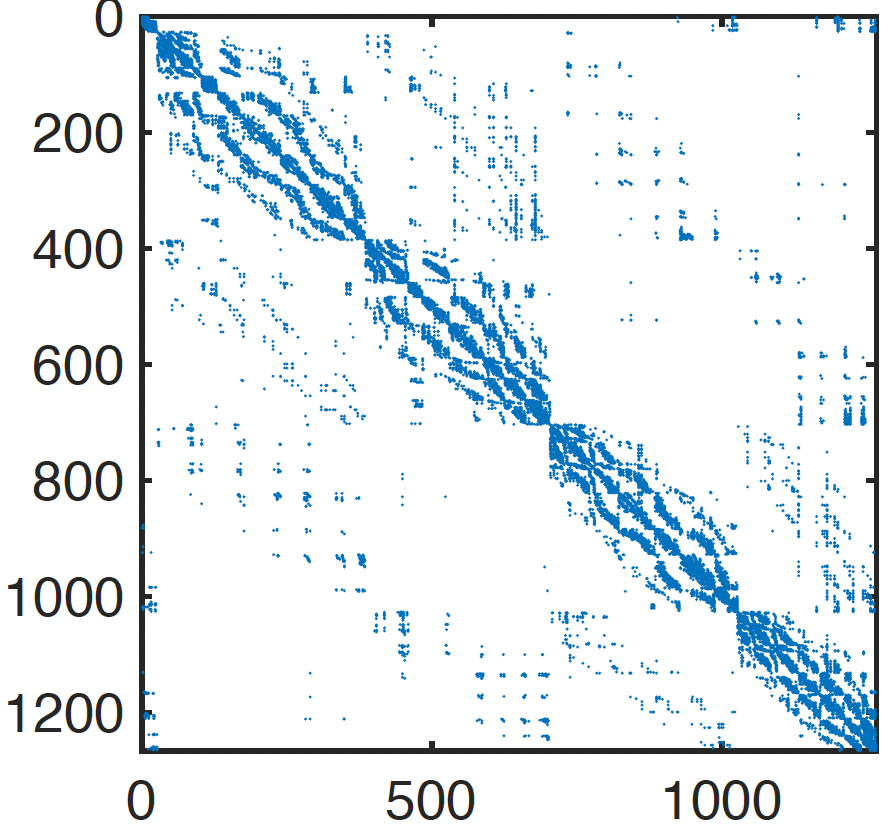}
\subcaption{Fourth Level.\\ $\kappa_2(U) =$ 6.54e16,\\ $\kappa_2(L) =$ 2.51e7.}
\label{fig:LU4_N14k_fill}
\end{subfigure}
\end{center}
\caption{Non-zero patterns of $L+U$ for matrix size $N=14186$ for the first
four levels using AMGToolbox. PeleLM model. Drop tolerance is set to 1.e$-15$ and 
fill limit per row set to 10.}
\label{fig:LU14k_fill}
\end{figure}

Figure \ref{fig:LU14k} displays the sparsity pattern of $L+U$ 
on the first four levels of AMG
using the AMGToolbox \cite{Verbeek2002} for matrix dimension $N=14186$.
Here, the drop tolerance is set to 1.e$-15$ and fill limit to $200$ per
row using the  ILUTP implementation in \cite{Carr2021} with pivoting turned 
off.\footnote{The term {\it pivoting} refers to row or column 
exchanges employed when a small {\it pivot}, or divisor, 
is encountered in the factorization.  With pivoting turned off, 
the pivot is always the diagonal element.} 
{A very small drop tolerance and large fill-in 
are unreasonable choices for an ILU factorization because they can 
substantially increase the cost associated with computing, applying and, storing 
the factors.  The plots emphasize the 
resulting vertical striping associated with poorly chosen parameters, which}
is attributed to small pivots combined with large amounts of fill-in. 
The solution is to 
enforce a smaller fill level per row \cite{ChowSaad1997}.  
For a very ill-conditioned $A$, limiting the fill level per row may be 
insufficient to prevent ill-conditioned $L$ and $U$. {However, 
enforcing a fill level that is too conservative may result in a poor 
approximation of $A$.}
Figure \ref{fig:LU14k_fill} displays the non-zero pattern of $L+U$ for the same 
matrices as in Figure \ref{fig:LU14k}, 
but with the fill level per row now set to $10$.  
The dramatic striping pattern is no longer present, however some remains, 
indicating the potential for ill-conditioning despite 
restricting the amount of fill.

{An example demonstrating the effects of 
a conservative drop tolerance on the conditioning of $A$ is omitted here, 
however, by imposing a very small drop tolerance, the 
resulting triangular solves may be unstable.  This is due to the fact that
the off-diagonal elements are much larger than those along the
diagonal\cite{ChowSaad1997}.}  
In fact, the 
prescription given in \cite{ChowSaad1997} is to consider scaling to 
reduce $\kappa(A)$ with the warning that non-normal
triangular factors may result. 
{In Section \ref{sec:analysisLU}, we 
demonstrate the effects of varying drop tolerances, combined with scaling 
of a judicious choice in drop tolerance followed by row/column scaling of 
factors results in very few iterations required to accurately solve the 
triangular systems.}

\section{An ILU Smoother with Scaled Factors}\label{sec:ILUSmooth}
{To mitigate the high degree of non-normality of $U$, we 
propose row and column scaling, as well as row scaling alone, of the triangular 
factors within the smoothing step of AMG.  We note again that in the applications 
considered, the $L$ factor generally did not display ill-conditioning 
or high degrees of non-normality, and so we focus our discussion on $U$. By scaling 
$U$, and subsequently reducing both $\kappa(U)$ and dep$(U)$, we facilitate the use 
of the significantly faster SpMV products appearing in the Richardson iteration on 
GPU architectures to solve the triangular linear systems. We first introduce the 
general form of the scaled linear system and demonstrate the effects on $\kappa(U)$ 
and dep$(U)$ for several applications; see Section \ref{sec:Scaling}. Then, we expand 
the Richardson iteration explicitly to derive the Neumann series, demonstrating 
how the series can be truncated to very few terms; see Section \ref{sec:analysisLU}.} 

\subsection{{Equilibration} to Reduce Non-Normality}

\label{sec:Scaling}
Row and column scaling is a common choice for reducing ill-conditioning, 
but also extends to non-normal matrices arising in the ILU factorization of a highly 
ill-conditioned matrix.  In doing so, we can expect the weight of the off-diagonal elements 
to decrease, particularly in the case of diagonally dominant matrices. 
Diagonal matrices $D_l$ and $D_r$ typically employ row or column norms.
Specifically, row and column scaling results in a linear system that has been 
equilibrated and now takes the form
\begin{equation}\label{eq:rowcolScale}
L\:D_r\:U\:D_c \:x = b.
\end{equation}
{The Richardson} iterations are applied to the upper 
triangular matrix $D_r\:U\:D_c$ (in addition to the lower triangular factor $L$). 
Because $D_r\:U\:D_c$ has a unit diagonal,  the iterations 
can be expressed in terms of a Neumann series, which is discussed 
{next, in Section \ref{sec:analysisLU}}. 
When an $LDU$ or $LDL^T$ factorization is available, the diagonal matrix $D$  
can represent row scaling for either the $L$ or $U$ matrix.  
For the applications considered, only the $U$ matrix is scaled, 
and thus row scaling is employed. Given an incomplete LU factorization, 
$D$ can be written as
\begin{equation}\label{eq:extractD} 
D = diag(U),
\end{equation}
and the scaled $\widetilde{U}$ is subsequently defined as
\begin{equation}\label{eq:newD}
    \widetilde{U} = D^{-1}U
\end{equation}
to obtain the incomplete $LD\widetilde{U}$ factorization from the incomplete LU factorization.

{In Tables \ref{table:dep} 
and \ref{table:cond}, the departure from normality and 
condition number of the factors, respectively, are given. 
We also compare $\kappa(U)$ and dep$(U)$ before and 
after scaling for matrices obtained from the
SuiteSparse Matrix Collection \cite{Davis2011} (first five rows).
Results for three matrices exported from PeleLM 
(final three rows) are also provided. Note that $\text{dep}(L)$ 
is generally very modest, and never as large as $\text{dep}(U)$. To generate these results
Matlab's {\tt ilu} was applied with setup type `nofill' (i.e.
ILU$(0)$).}

\begin{table}[h!]
\centering
\begin{tabular}{|l|c|c|c|c|c|}
\hline
\textbf{Matrix}&Dimension&$\text{dep}(L)$&$\text{dep}(U)$&$\text{dep}(D^{-1}U)$&$\text{dep}(D_rUD_c)$\\
&&&&(row scaling)&(row/col scaling)\\\hline
{\tt af\_0\_0\_k101} & 503625 & 326.95 & 1.84e8 & 326.95 & 320.89 \\\hline
{\tt af\_shell1} & 504855 & 386.66 & 1.52e8 & 386.66 & 407.35 \\\hline
{\tt bundle\_adj} & 513351 & 8.52e6 & 4.52e11 & 8.52e6 & 438.70 \\\hline
{\tt F1} & 343791 & 335.52 & 4.89e8 & 335.52 & 331.79 \\\hline
{\tt offshore} & 259789 & 231.86 & 7.05e15 & 231.86 & 222.71 \\\hline 
PeleLM331 & 331 & 8.37 & 1.50e6 & 8.37 & 4.13 \\\hline 
PeleLM2110 & 2110 & 16.99 & 1.09e7 & 16.99 & 9.33 \\\hline 
PeleLM14186 & 14186 & 1.45e4 & 1.00e6 & 1.45e4 & 26.33 \\\hline
\end{tabular}
\caption{Departure form normality for the $L$ and $U$ factors when applying an ILU(0) factorization to several matrices, followed by the departure from normality after row scaling and row/col scaling are applied to the $U$ factor. The first five come from \cite{Davis2011}, and the last three are extracted from PeleLM.  Matlab's {\tt ilu} with type `nofill' was computed for all matrices.}
\label{table:dep}
\end{table}

\begin{table}[h!]
\centering
\begin{tabular}{|l|c|c|c|c|c|}
\hline
\textbf{Matrix}&$\kappa(A)$&$\kappa(L)$&$\kappa(U)$&$\kappa(D^{-1}U)$&$\kappa(D_rUD_c)$\\
&&&&(row scaling)&(Ruiz scaling)\\\hline
{\tt af\_0\_0\_k101} & 3.60e8 & 156.54 &1.02e3  & 75.78 & 108.83 \\\hline
{\tt af\_shell1} & 1.72e10 & 49.99 & 231.94 & 116.42 & 171.30 \\\hline
{\tt bundle\_adj} & 6.10e15 & 4.59e12 & 3.53e14 & 2.93e12 & 2.37e3  \\\hline
{\tt F1} & 3.26e7 & 6.51e3 & 1.34e5 & 2.67e4 & 1.54e4 \\\hline
{\tt offshore} & 2.32e13 & 96.79 & 7.56e10 & 148.35 & 156.99 \\\hline
PeleLM331 & 3.48e17 & 14.06 & 3.87e9 & 43.06 & 18.72 \\\hline
PeleLM2110 & 3.21e17 & 13.39 & 4.41e9 & 34.02 & 12.31 \\\hline
PeleLM14186 & 6.64e15 & 1.83e8 & 6.87e12 & 1.74e7 & 9.51 \\\hline
\end{tabular}
\caption{{Condition number for $A$, $L$, and $U$ when applying an ILU(0) factorization to several matrices, followed by the condition numbers after row scaling and Ruiz scaling are applied to the $U$ factor. The first five come from \cite{Davis2011}, and the last three are extracted from PeleLM.  Matlab's {\tt ilu} function with type `nofill' was computed for all matrices.}}
\label{table:cond}
\end{table}

Our results demonstrate that the 
scaling strategies substantially reduce $\text{dep}(U)$, and that in many cases, the 
row/column scaling produces a larger 
reduction compared with row scaling alone. 
However, in some cases, 
row scaling results in a lower departure from normality.  
When row scaling results in a smaller
$\text{dep}(U)$, both are nearly the same, or at least of the same order
of magnitude. 
The ILU smoother with scaling 
is provided in Algorithm \ref{alg:ILUJac}.
To employ row scaling within Algorithm \ref{alg:ILUJac}, 
formation of $D_U$, $D$ is constructed as in (\ref{eq:extractD}) 
and updated $\widetilde{U}$ as in (\ref{eq:newD}).  Here, only the $U$ factor
is scaled. The ILU factorization is assumed to result 
in a lower triangular matrix with unit diagonal, e.g. $L = I + L_s$
\footnote{This is a reasonable assumption as many factorizations 
produce such a matrix; e.g. consider the {\tt ilu} function in Matlab.}, 
which also generates a finite Neumann sum. {For the ILUT Schur complement smoother, we again mitigate high degrees of non-normality by limiting fill-in, and it is important to note the residual vector
is not required for the Schur complement GMRES 
solver for a fixed number of iterations.
The initial guess is set to $x^{(0)} = 0$ and $r^{(l)} = b$,
without a convergence check with the residual $r^{(1)}$.}

\begin{algorithm}
\centering
\begin{algorithmic}[0]
    \State{Given $A\in\mathbb{C}^{n\times n}$, $b\in \mathbb{C}^n$}
    \State{Define $droptol$ and $fill$}
  \State{Compute $A \approx LU$ with $droptol$ and $fill$ imposed}
  \State{Define $m_L$ and $m_u$, total number of iterations for solving $L$ and $U$}
  \State{Define $y = 0$, $v = y$}
  \State{//Richardson iteration to solve $Ly = b$}
  \For{$k = 1:m_{L}$}
    \State{$y = b - L_s\:y$ }
  \EndFor 
  \State{Scale  $U$ and $y$ to obtain 
  scaled $\widetilde{U}$ and $\widetilde{y}$, and $D_k$}
  \State{Let $D_U = diag(\widetilde{U})$}
  \State{Define $D = D_U^{-1}$}
  \State{//Richardson iteration to solve $\widetilde{U}v = \widetilde{y}$}
  \For{$k = 1:m_U$}
    \State{$v = D\:\widetilde{y} - \widetilde{U}_s\:v$ }
  \EndFor 
  \State{//Update and unpermute the solution}
  \State{$v = Dv$}
  \State{$x = P^{-1}v$}
\end{algorithmic}
\caption{\label{alg:ILUJac} ILU+Richardson smoother for 
AMG with row scaling of $U$.
}
\end{algorithm}

\subsection{Neumann Series and the Richardson Iteration}\label{sec:analysisLU}
{Our technique for} scaling the triangular factors not only reduces the 
departure from normality, but also results in 
a finite Neumann series. 
{First consider the Jacobi iteration for solving the linear system 
$Ax = b$.  Note that Jacobi has been discussed previously as a smoother, but it is a 
well known iterative solver for an arbitrary linear system (unrelated to AMG).  
Given the preconditioner  $M = D$, the non-compact form of (\ref{eq:compact}) is given by}
\begin{equation}
x^{(k+1)} =x^{(k)} + D^{-1}\:(\:b - A\:x^{(k)}\:)
\end{equation}
{where $D$ is the diagonal part of $A$. For the triangular systems
resulting from the ILU factorization $A \approx LDU$ (as opposed to
the regular splitting $A = D + L + U$), the iteration matrices 
are denoted $G_L$ and $G_U$ for the lower and upper triangular
factors, $L$ and $U$, respectively.  Let $D_L$ and $D_U$ be the diagonal parts of 
the $L$ and $U$ and let $I$ denote the identity matrix.
Assume $L$ has a unit diagonal, then}
\begin{eqnarray}
G_L & = & D_L^{-1}\:(\: D_L - L \:) = I - L, \\
G_U & = & D_U^{-1}\:(\:D_U - U\:) = I - D_U^{-1}\:U.\label{eq:LUiterMat}
\end{eqnarray}
{We next consider the effect of iterating 
with a scaled ILU factorization 
(e.g., ILU scaled with row/column scaling {as in 
(\ref{eq:rowcolScale})}, or the row-scaled LDU {using 
(\ref{eq:newD})}).} The iteration matrix (\ref{eq:LUiterMat})  simplifies to 
$G_U = U_s$, where $U_s$ is a strictly upper triangular matrix.
\footnote{{Again, recall we omit discussion of $L$, and thus 
$G_L$ for simplicity, as the applications demonstrate, do not exhibit 
problematic dep$(L)$.}} To solve $Ux = b$ let $b_s = D_U^{-1}b$, and 
$U = I + U_s$. Then replace (\ref{eq:compact}) with a Richardson iteration
\begin{equation}
x^{(k+1)} = b_s + (\:I - U\:)\: x^{(k)} = b_s - U_s\:x^{(k)}
\label{eq:richard}
\end{equation}
where the unit diagonal is removed. After expanding it follows that
\begin{eqnarray}
x^{(k+1)} & = & b_s - U_s\:b_s + U_s^2\:b_s - \dots + (-1)^k\:U_s^k\:b_s
    \nonumber \\
    & = & (\:I - U_s + U_s^2 - \dots + (-1)^k\:U_s^k\:)\:b_s
        \label{eq:Neumann} \\
    & = & (\:I+U_s\:)^{-1}\:b_s. \nonumber
\end{eqnarray}
Then, the inverse of $U$ is expressed as a Neumann series
\begin{equation}
U^{-1} = (\: I + U_s\:)^{-1} = I - U_s + U_s^2 - \cdots
= \sum_{i=0}^{n}(-1)^{i}U_s^{i},
\label{eq:Neuman}
\end{equation}
and with $U_s$ strictly upper triangular and nilpotent, 
the above sum is necessarily finite.  

The series in (\ref{eq:Neuman}) converges when 
$\|U_s\|_2 < 1$, and in practice,
this is true for the ILU$(0)$ and ILUT for certain drop 
tolerances, where a convergent Neumann series is guaranteed.  
Even in cases when $\|U_s\|_2 \geq1$, it is observed that $\|U_s^p\|_2 <1$ 
for $p\ll n$, permitting truncation of the Neumann series - 
and thus the Richardson iteration - to a small number of terms. 
That $\|U_s^p\|_2$ eventually decreases is not unexpected 
because the number of 
possible non-zeros in $U_s^p$ necessarily grows smaller as $p$ grows larger 
when $U_s$ is dense.  However, numerical nilpotence is observed 
for $p\ll n$ in many cases for sparse $U_s$, and the size of $p$ clearly 
depends on the number of non-zeros allowed in $U$ and consequently $U_s$ 
(either by imposition of small $droptol$ or conservative fill, or both).  
In other words, $\|U_s^p\|_2$ is effectively zero
much sooner than the theoretically guaranteed $\|U_s^n\|_2$. 
    
\begin{figure}[hh]
\captionsetup[subfigure]{font=large}
\begin{center} 
\begin{subfigure}{.49\textwidth}
\includegraphics[width=\linewidth]{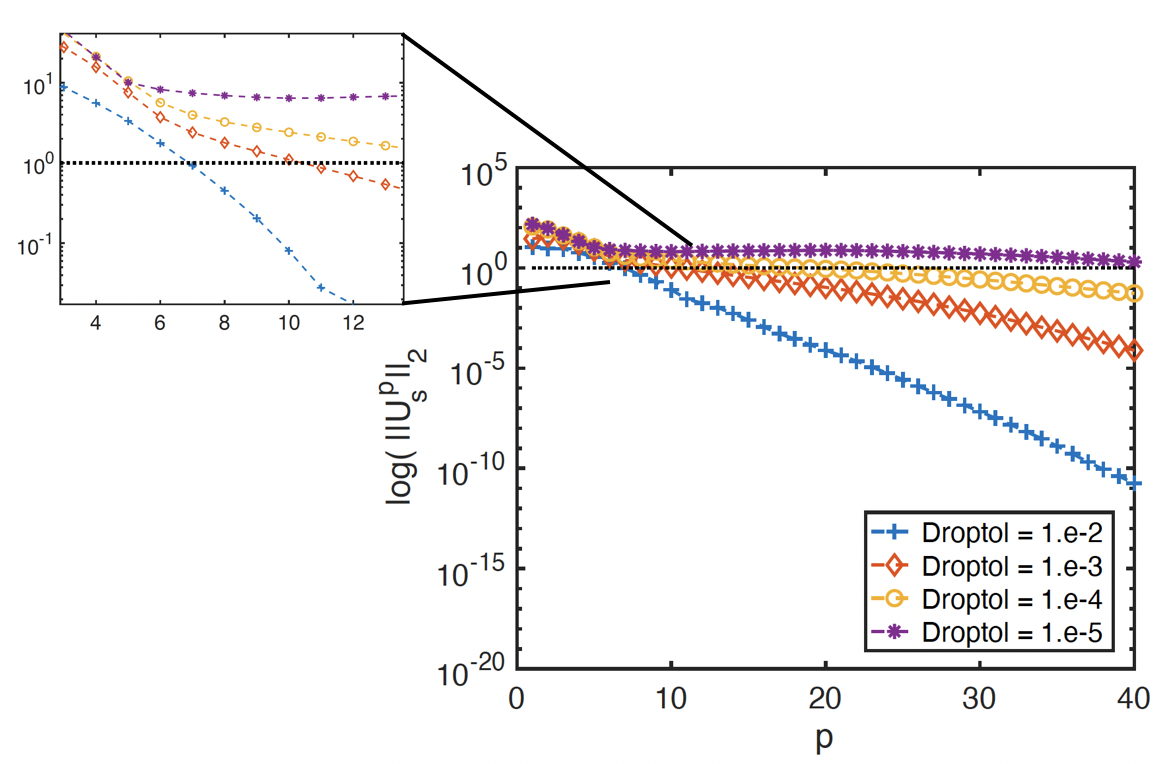}
\caption{{Using row scaling only.}} \label{fig:normUsLDU}
\end{subfigure}
\begin{subfigure}{.48\textwidth}
\includegraphics[width=\linewidth]{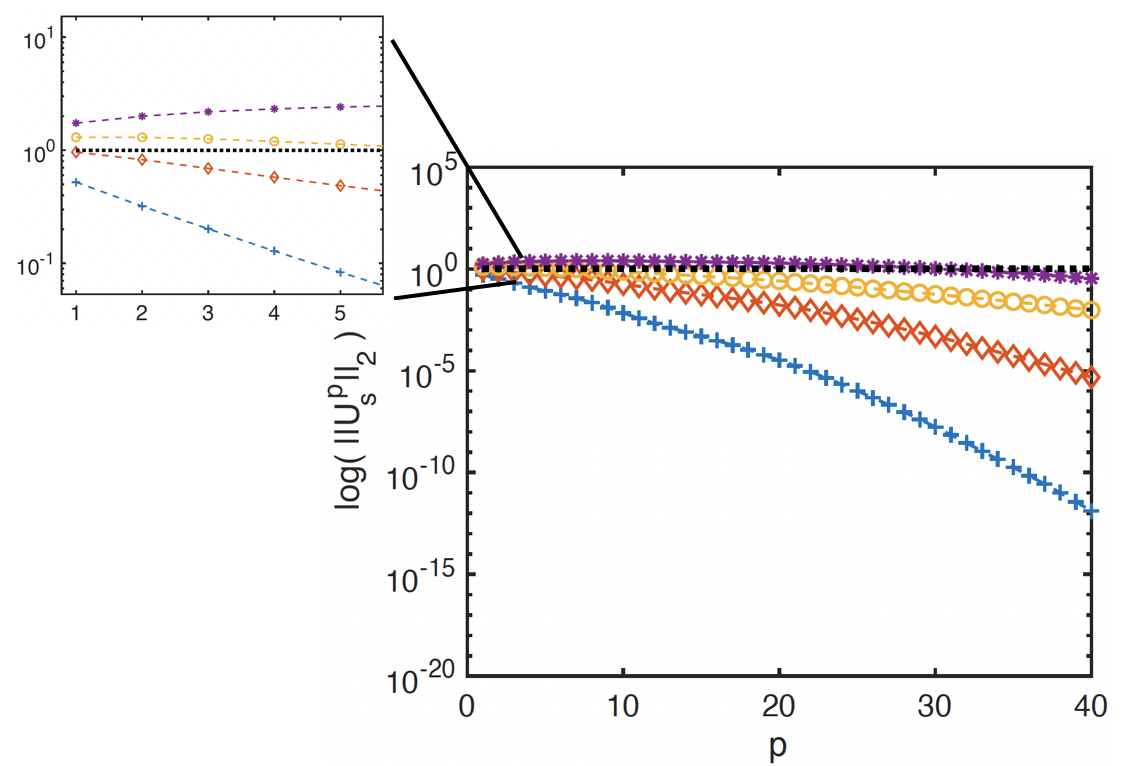}
\subcaption{{Using row and column scaling.}}
\label{fig:normUs}
\end{subfigure}
\end{center}
\caption{{$\|U_s^p\|_2$ for $p = 1, 2, \dots 40$, for $U = I + U_s$ using scaling and where $U_s$ is the strictly upper triangular part of $U$. 
The black, dotted line represents the bound $1$. PeleLM matrix dimension $N=14186$.}}
\label{fig:Up}
\end{figure}

{Figure \ref{fig:Up} displays $\|U_s^p\|_2$ when applying 
scaling to the ILU smoother employed at the finest level. Here, Matlab's 
{\tt ilu} is used with type `ilutp', 
threshold $0$ (i.e. no pivoting), and various drop tolerances. 
For larger drop tolerances (i.e. $droptol = 1.e-2$ for both row/column 
and row scaling), $\|U_s^p\| < 1$ for $p = 1$, giving convergence of the 
Neumann series.  However, not for smaller drop tolerances. 
$\|U_s^p\|_2$ actually increases for the
first few values of $p$, but eventually decreases and falls below $1$. 
Figure \ref{fig:normUsLDU}, displays $\|U_s^p\|_2$ for matrix dimension 
$N = 14186$ with row scaling. $\|U_s^p\| < 1$ 
for modest $p$ (e.g. $p = 7$ for $droptol = 1.e-2$).  For smaller drop
tolerances, $p$ can be moderately large (e.g. $p=45$ for $droptol = 1.e-2$).}

The strategy is now extended to much larger linear systems 
exported from the Nalu-Wind \cite{SC21}, and PeleLM \cite{PeleLM} 
pressure solvers and a parallel performance analysis is provided. 
The Nalu-Wind matrix has dimension $N=21$M and is
sufficiently large to exhibit differences in the strong scaling
characteristics of the ILU Schur complement smoother. The PeleLM
system of dimension $N = 11$M is highly ill-conditioned and requires
${\cal O}(100)$ GMRES+AMG iterations to converge with a 
Gauss-Seidel smoother.

\section{Numerical Results}\label{sec:results}

Our iterative approach is applied to the solution of linear systems
from the PeleLM \cite{PeleLM} ``nodal projection'' step and the pressure 
continuity equation for the Nalu-Wind \cite{SC21} CFD model. 
Relatively large matrices of dimension, greater than 10M
were exported from the PeleLM nodal pressure projection
solver \cite{PeleLM} and the Nalu-Wind pressure continuity solver \cite{SC21}. 
Furthermore,  scaling was applied directly to the non-normal $U$ factor to 
reduce its departure from normality, facilitating the use of a Neumann series
(Richardson iteration) to solve the triangular system.  A performance model is also
provided together with numerical and parallel scaling results 
in the following sections.

The stopping criteria for Krylov methods is an important
consideration and is related to backward error for
solving linear systems $Ax = b$. The most
common convergence criterion found in existing iterative
solver frameworks is based upon the relative residual,
defined by
\begin{equation}
\frac{\|r_k\|_2}{\|b\|_2} =
\frac{\|b - Ax_k\|_2}{\|b\|_2} < tol,
\label{eq:relres}
\end{equation}
where $r_k$ and $x_k$ represent, respectively, the residual and 
approximate solution after $k$ iterations of the iterative solver.
An alternative metric commonly employed in direct solvers
is the norm-wise relative backward error (NRBE)
\begin{equation}
\textit{NRBE} = \frac{\|r_k\|_2}{\|b\|_2 + \|A\|_{2}\|x_k\|_2}.
\label{eq:nbre}
\end{equation}
In numerical experiments, the norm-wise relative backward error 
for the solution of linear systems with GMRES was sometimes found to
be lower than when the right-preconditioned GMRES was employed \cite{Saad86}.
Indeed, the latter exhibited false convergence (the Arnoldi residual
norm did not agree with the norm of the true residual 
$r^{(k)} = b - Ax^{(k)}$) when
executed in parallel for highly ill-conditioned problems, $\kappa(A) > 10^{15}$.
Flexible FGMRES \cite{Saad03} was found to be the most effective Krylov solver 
with an AMG preconditioner and did not exhibit false convergence.

The convergence of GMRES+AMG using ILU smoothers is compared with polynomial Gauss-Seidel.  
In particular, ILU on the finest level is combined with Gauss-Seidel 
on coarser levels. The choice to apply the ILU smoother on 
any number of levels -- starting from the finest level -- is now 
an option available  in hypre \cite{Falgout2002}. 
The ILU smoother on all levels and the ILUT Schur 
complement smoother are evaluated.
{The hypre-BoomerAMG
library was designed for massively-parallel computation 
\cite{Baker:2011} and now also supports GPU acceleration 
of key solver components \cite{Falgout2021}.
The new formulation allows simple and efficient implementations that can utilize  available optimized sparse kernels on GPUs \cite{SC21}.
Li et. al.~\cite{Li2020} describe the class of M-M based interpolation operators suitable for efficient GPU computations.}

\subsection{PeleLM Combustion Model}\label{sec:peleresults}

\begin{figure}[hbt!]
\centering
\includegraphics[height=0.35\textheight,width=0.65\textwidth]{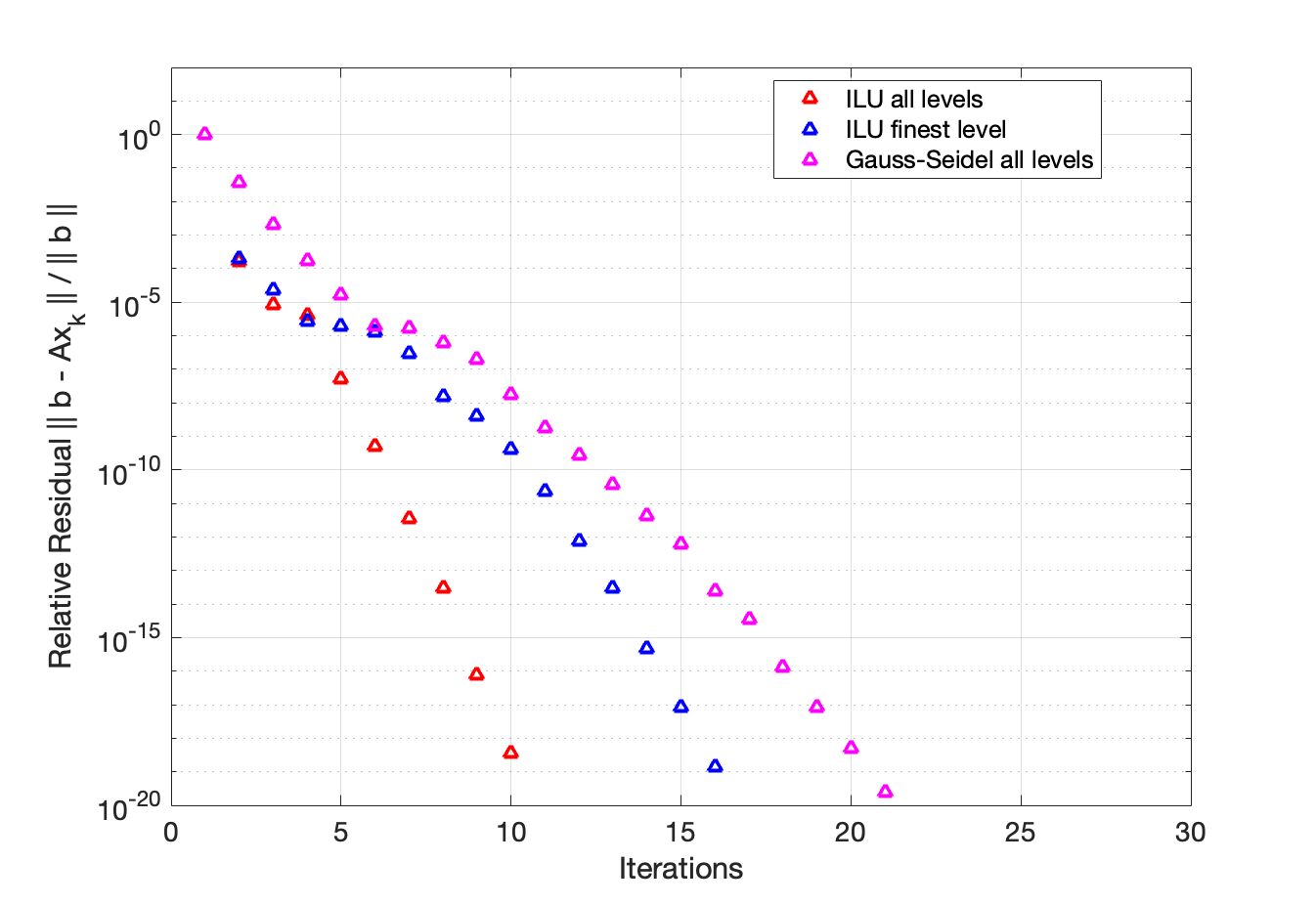}
\caption{\label{fig:multi-level-ilu}  
AMGToolBox results. PeleLM model.
ILUT smoother on finest level versus all levels.
Matrix size $N=14186$.}
\end{figure}

A sequence of three different size problems was examined, based on matrices
exported from the PeleLM pressure continuity solver \cite{PeleLM}. 
The problem being solved is combustion in a piston-cylinder configuration,
where the curved surface of the cylinder requires cut-cells through the mesh.
The first of these matrices is a dimension $N=14186$ linear system, solved with
GMRES$+$AMG using the AMGToolBox \cite{Verbeek2002, Joubert2006} by applying ILUT
smoothing only on the finest level $\ell=1$, then ILUT on all levels
and polynomial Gauss-Seidel smoothers. Iterative Richardson 
solvers are employed. Two pre- and post-smoothing 
sweeps were applied on all $V$-cycle levels, except for the coarse
level direct solve. The AMG strength of connection 
threshold was set to $\theta=0.25$. The convergence histories
are plotted in Figure \ref{fig:multi-level-ilu}. The lowest
iteration count results from using ILU(0) smoothing on all
levels, however, the minimum compute time is obtained by
using ILU(0) only on the finest level and polynomial Gauss-Seidel
on the remaining levels. The iterative triangular solves fail
to converge unless either a preconditioned Jacobi or Richardson
iteration is employed. In the case of ILU(0), 
both row and row/column scaling exhibit similar convergence
histories and thus row scaling is less costly in the set-up phase.

In the case of ILUT as the smoother, there are differences in the achievable GMRES
error level between the row and row/column scaling, depending on the drop tolerances
and level of fill per row. For the dimension $N=14$K problem,
the norm-wise relative backward error (NRBE) is reported in 
Table \ref{tab:NRBE2} for $\textit{droptol} = 1$e$-3$ and $\textit{lfill} = 5$,
where the NRBE is found to be an order of magnitude lower.
A slightly lower backward error was obtained using two Richardson
iterations with the row/column scaling versus three iterations using row scaling
may not be sufficient to justify the additional set-up cost.


\begin{table}[htb]
\centering
\begin{tabular}{|l|c|c|c|c|c|} \hline
{\bf Richardson Iters.}  &    5 &     4  &    3 &      2 &     1 \\ \hline
{\bf Row scale $U$}      &  5.85e$-10$ & 7.3e$-10$ & 5.94e$-10$ & 3.84e$-10$ & 5.35e$-8$ \\
\hline
{\bf Row/col scale $U$}  &  5.85e$-10$ & 7.17e$-10$ & 5.93e$-10$ & 3.85e$-10$ & 1.91e$-9$\\ 
\hline
\end{tabular}
\caption{ \label{tab:NRBE2} 
GMRES$+$AMG norm-wise relative
backward error (NRBE) when using ILU+Richardson smoother with a varying number of Richardson iterations for the PeleLM matrices with dimension $N=14186$, $\textit{droptol} = 1$e$-3$, and
$\textit{lfill} = 5$. Row and Row/column scaling is applied to the $U$ factor only.}
\end{table}

\begin{figure}[htb!]
\centering
\includegraphics[height=0.35\textheight,width=0.65\textwidth]{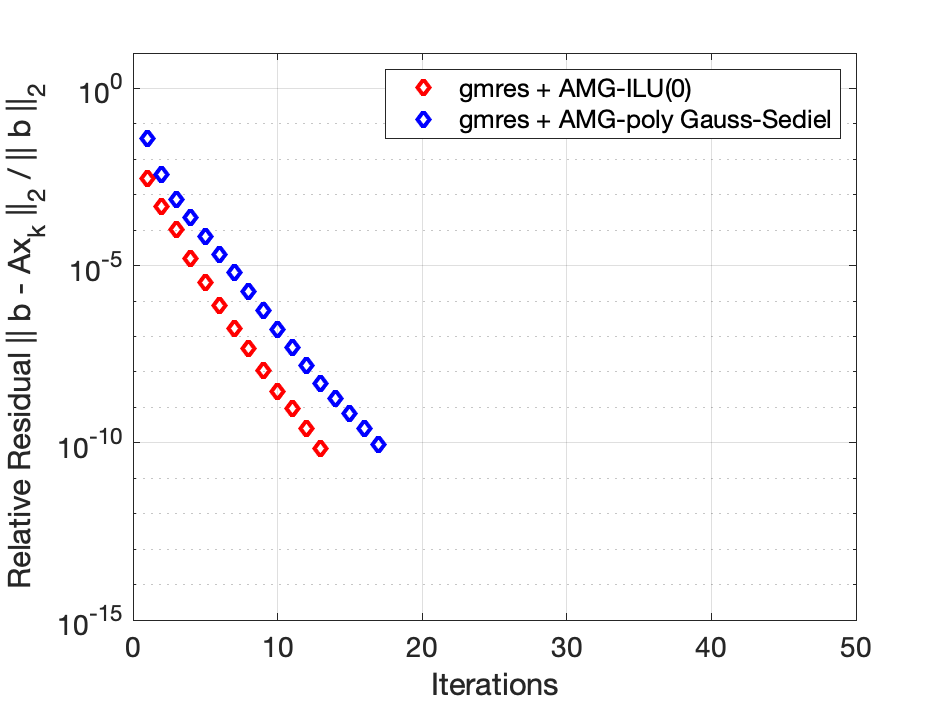}
\caption{\label{fig:1.4million} 
hypre-BoomerAMG GPU results  NREL Eagle.
Convergence history of (F)GMRES+AMG 
with polynomial, and  mixed ILU smoothers
with iterative solves. PeleLM model. 
Matrix size $N=1.4$M}
\end{figure}

Results using hypre-BoomerAMG for the $N=1.4$M 
linear system  are plotted in Figure \ref{fig:1.4million}.
These tests were performed on the NREL Eagle supercomputer with 
Intel Skylake CPUs and NVIDIA V100 GPUs. The parallel maximum independent set
(PMIS) algorithm is applied together with aggressive coarsening 
and ``MM-ext+i'' interpolation are employed, with a strength of
connection threshold, $\theta = 0.25$.  Because
the problem is very ill-conditioned, flexible FGMRES achieves the 
best convergence rates and the lowest NRBE. Iterative solvers
were employed in these tests with ten (10) iterations.
The convergence histories are plotted for mixed-ILUT,
and polynomial Gauss-Seidel smoothers. The ILUT parameters
were $\textit{droptol} = 1$e$-2$ and $\textit{lfill} = 10$. The lowest 
time for a  single-GPU, was the ILUT smoother 
with Richardson iterations which achieved a solve time of $0.11$ seconds.

\begin{figure}
\centering
\includegraphics[height=0.35\textheight,width=0.65\textwidth]{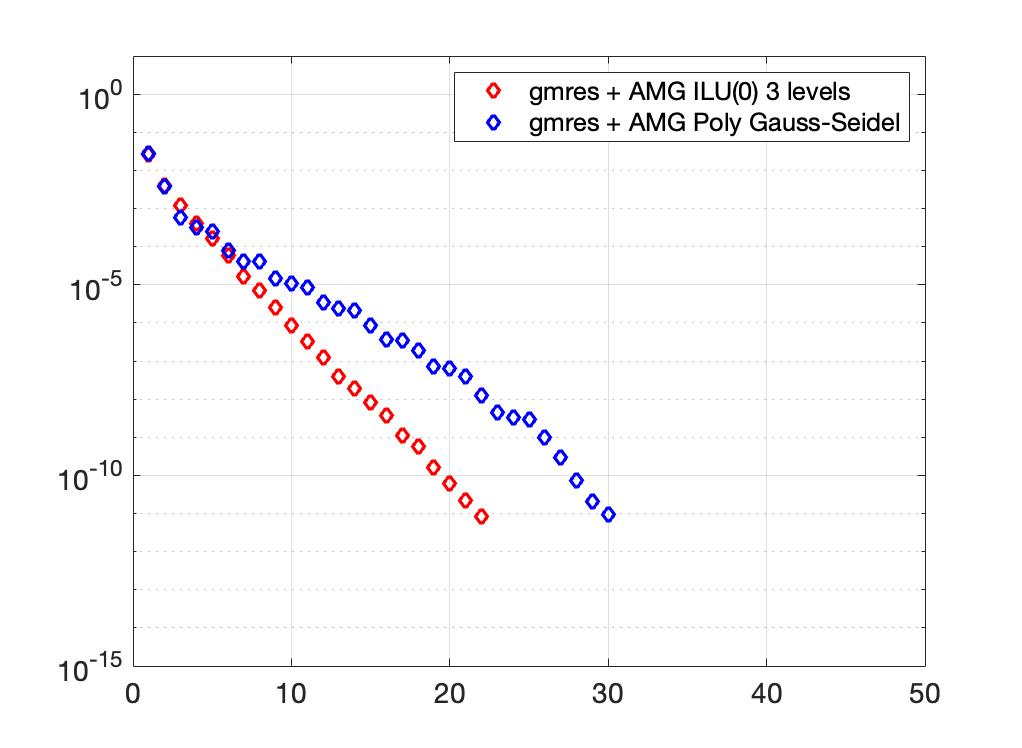}
\caption{\label{fig:4million} 
hypre-BoomerAMG GPU results  NREL Eagle.
Convergence history of (F)GMRES+AMG 
with polynomial, and  mixed ILU smoothers
with iterative solves. PeleLM model. Matrix size $N=4$M}
\end{figure}

For comparison, a larger PeleLM linear system of dimension $N= 4$M
was run on the ORNL Crusher supercomputer (see Figure \ref{fig:4million}). 
The machine contains multiple
compute nodes consisting of two AMD EPYC CPU sockets and four AMD MI250X 
GPU sockets, each containing two GCDs (GPU compute devices) for eight total.
The same AMG parameters as the dimension $1.4$M problem were specified, however,
and ILU($0)$ smoother is applied on the first three $V$-cycle levels. The 
Richardson iterations are reduced to six (6) upper and five (5) lower per level.
To reduce the relative residual to 1e$-11$, the solve time is 0.11 seconds
with iterative triangular solvers. Whereas, the solve time is 0.16 seconds
when a direct solver is employed in the smoother. The speed-up is now
reduced to $1.5\times$ on the AMD MI250X GPU and may be attributed to
a faster implementation of the direct triangular solver by AMD.
Despite the increased number of GMRES+AMG iterations for the polynomial
Gauss-Seidel smoother, the solve time is $0.16$ seconds, which is
comparable to the direct solves with ILU. However, the number of
GMRES iterations has grown from the smaller problem and is expected
to increase further at larger problem sizes.

\begin{figure}[htb!]
\centering
\includegraphics[height=0.35\textheight,width=0.65\textwidth]{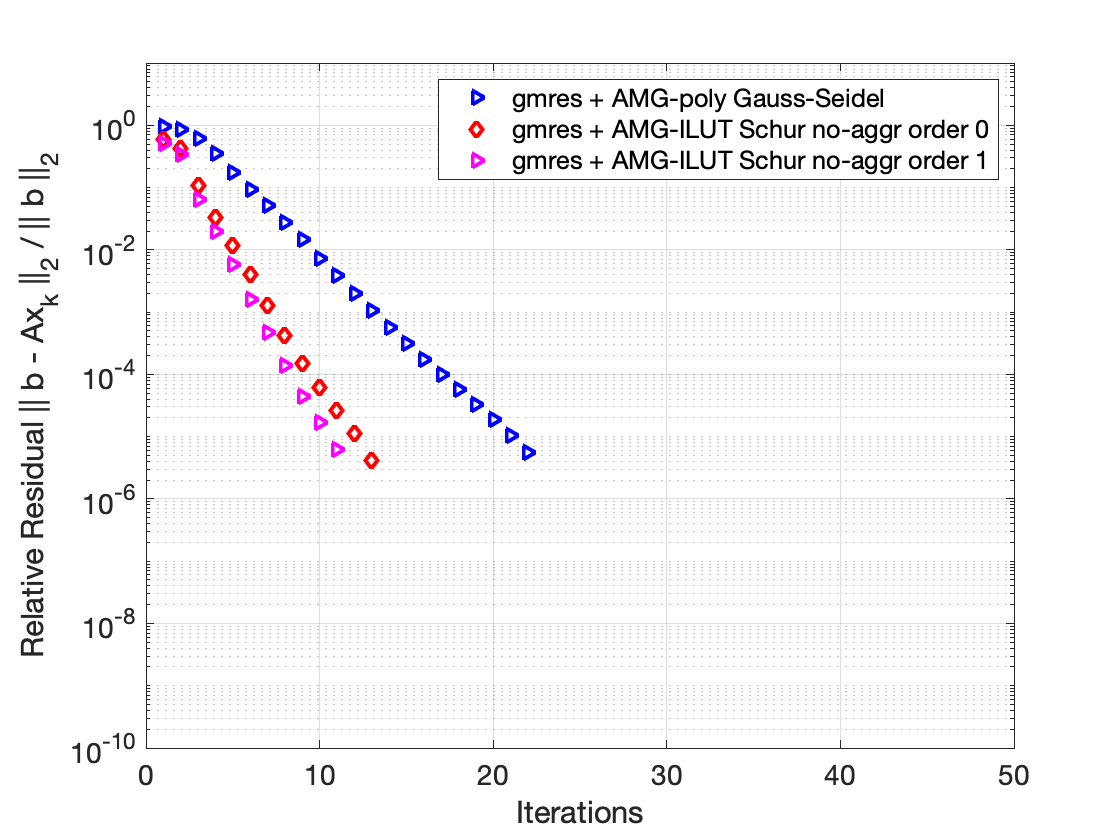}
\caption{\label{fig:23million} hypre-BoomerAMG GPU results
an NREL Eagle. Nalu-Wind.
Convergence history of (F)GMRES+AMG 
with polynomial, and ILUT Schur Complement smoothers
with iterative solves.
Matrix size $N=23$M}
\end{figure}

\subsection{Exa-Wind Fluid Mechanics Models}

In Nalu-wind, the pressure systems are solved using MGS-GMRES with an
AMG preconditioner, where a polynomial Gauss-Seidel
smoother is described in Mullowney 
et al.~\cite{SC21}.  Hence, Gauss-Seidel is a compute time intensive
component, when employed as a smoother within an AMG $V$-cycle.
The McAlister experiment for wind-turbine blades is an
unsteady RANS simulation of a fixed-wing, with a NACA0015 cross section, 
operating in uniform inflow. 
Resolving the high-Reynolds number boundary layer over the wing surface 
requires resolutions of ${\cal O}(10^{-5})$ normal to the surface resulting in 
grid cell aspect ratios of ${\cal O}$(40k). These high aspect ratios present a 
significant challenge. Overset meshes were employed
to generate body-fitted meshes for the wing and the wind tunnel geometry.
The simulations were performed for the wing at a 12 degree angle of attack, 1m chord length, 
denoted $c$, 3.3 aspect ratio, i.e., $s = 3.3c$, and square wing tip. 
The inflow velocity is $u_{\infty}= 46$ m/s, the density is $\rho_{\infty} = 1.225$
${\rm kg/m}^3$, and dynamic viscosity is $\mu = 3.756\times 10^{-5}$ kg/(m s), 
leading to a Reynolds number, $Re = 1.5\times 10^6$.
Wall normal resolutions were chosen to adequately represent the boundary layers 
on both the wing and tunnel walls. The $k - \omega$ SST RANS turbulence model 
was employed for the simulations. Due to the complexity of mesh generation, 
only one mesh with approximately 3M grid points was generated. 

\begin{figure}[htb!]
\centering
\includegraphics[width=0.55\textwidth]{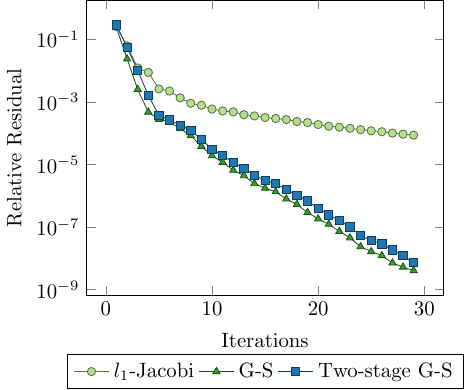}
\caption{\label{fig:conmat} GMRES$+$AMG with $\ell_1$--Jacobi, Gauss-Seidel 
and polynomial Gauss-Seidel smoothers for a linear system from Nalu-Wind model with dimension $N=3.1$M.}
\end{figure}

Coarsening is based on the parallel maximal
independent set (PMIS) algorithm 
allowing for a parallel setup phase.  
The strength of connection threshold is set to $\theta = 0.25$.  
Aggressive coarsening is applied on the
first two $V$-cycle levels with multi-pass 
interpolation and a stencil width of
two elements per row.  
The remaining levels employ M-M extended$+$i interpolation,
with truncation level $0.25$ together with a maximum
stencil width of two matrix elements per row. 
The smoother is hybrid block-Jacobi with two sweeps of polynomial
Gauss-Seidel applied locally on an MPI rank and then Jacobi 
smoothing for globally shared
degrees of freedom.  The coarsening rate for the wing simulation is roughly
$4\times$ with eight levels in the $V$-cycle for hypre.  Operator complexity
$C$ is close to $1.6$ indicating more efficient $V$-cycles with aggressive
coarsening, however, an increased number of GMRES iterations 
are required compared to standard coarsening. 
The comparison among $\ell_1$--Jacobi,  Gauss-Seidel and 
the polynomial Gauss-Seidel smoothers is shown in Figure~\ref{fig:conmat}.

In order to evaluate the differences between the direct and iterative 
triangular solvers employed in the smoother,
the compute times for a single GMRES$+$AMG pressure solve
are given in Table \ref{tab:smoothers}. The $\ell_1$--Jacobi
smoother from hypre is included for comparison. Both the
CPU and GPU times are reported for the NREL Eagle
supercomputer with Intel Skylake Xeon CPUs and
NVIDIA V100 GPU's. In all cases, one sweep of
Gauss-Seidel and two sweeps of $\ell_1$--Jacobi are employed
because the number of sparse matrix-vector multiplies
are equivalent in both cases. Either one CPU core or GPU 
was employed in these tests. The time reported corresponds
to when the relative residual has been reduced below 1$e-$5.

\begin{table}[htb]
\centering
\begin{tabular}{|l|c|c|c|} \hline
 & $\ell_1$--Jacobi & Gauss-Seidel & Poly G-S \\ \hline
iterations &   36 &  18 &   21 \\ \hline
CPU (sec)  & 19.4 &  12 & 19.5 \\ \hline
GPU (sec)  & 0.37 & 3.2 & 0.27 \\ \hline
\end{tabular}
\caption{ \label{tab:smoothers} GMRES$+$AMG compute time. 
Jacobi, Gauss-Seidel and Polynomial Gauss-Seidel Smoothers for a linear system from Nalu-Wind model with dimension $N = 3.1$M.}
\end{table}

The timing results indicate the solver time with Gauss-Seidel
is lower than when the $\ell_1$ Jacobi smoother is employed on
the CPU. However, the latter is more computationally efficient 
on the GPU. Whereas the polynomial Gauss-Seidel smoother leads 
to the lowest compute times on the GPU and results in a 
10X speed-up compared to the smoother
that employs a direct triangular solver.

The NREL 5-MW turbine is a notional reference 
turbine with a 126 meter rotor that is appropriate 
for offshore wind studies. 
These models use inflow and outflow boundary 
conditions in the directions
normal to the blade rotation and symmetry boundary conditions in other 
directions. For each simulation, 
50 time steps are taken from a cold start 
with four Picard iterations per time step.  
Convergence histories for one such pressure 
linear system (after reaching steady-state) are displayed in 
Figure \ref{fig:23million} for the ILUT-Schur complement 
and polynomial Gauss-Seidel smoothers. The former requires half
as many iterations to reach the 1e$-5$ convergence tolerance.
Furthermore, a coarse-fine (C-F) ordering of the degrees of
freedom results in fewer iterations.  The strength of connection
parameter was set to $\theta = 0.57$, which contributes to a
reduction in the AMG set-up time. In addition, the ILU drop
tolerance was $\textit{droptol} = 1.0 \times 10^{-2}$, 
with a fill level per row of $\textit{lfill} = 2$.
Sufficient smoothing was achieved with 18 iterations for the
lower triangular $L$ solve and 31 for the $U$ solve.
Aggressive coarsening was not specified and a single level
of ILU smoothing was applied.

\subsection{{$V$-Cycle GPU Performance Model}}

A mixed AMG algorithm is obtained with ILU smoothing
on the finest levels of the AMG $V$-cycle hierarchy (e.g. level 1),
followed by polynomial Gauss-Seidel {or $\ell_1$-Jacobi} on the
remaining levels. In the numerical results reported earlier, the
ILU smoothing is used on all levels and also in the mixed
configuration for comparison. The latter requires fewer sparse
matrix-vector multiplies and thus is more efficient.

{The floating point operation count} of the $V$-cycle, 
is determined by the number of non-zeros in the factors.
The ILU smoother requires twenty (20) Richardson iterations and one outer sweep.
These are applied during {\it pre-} and {\it post-smoothing} or twice. Therefore,
the total number of flops required for the $V$-cycle with ILU$(0)$ smoothing 
on every level is given by the sum below, where an SpMV costs 
$2\times \text{nnz}(A_l)$ on level $l$, 
\begin{equation}
    \text{flops} = \sum_{l=1}^{N_l-1}\: \text{nnz}(A_l) \: \times 80
\end{equation}
whereas the factor $80$ is replaced by $8$ for the polynomial smoother
of degree two. This makes a compelling case for the mixed AMG
approach with a combination of ILU and {$\ell_1$-Jacobi} smoothers
when the convergence rate is improved and leads to lower flop counts.
The cost of the coarse grid direct solve is ${\cal O}(N_c^3)$, where
$N_c$ is the dimension of the coarsest level matrix $A_c$, and is small
in comparison. The cost of the Krylov iteration is dominated by the
SpMV with the matrix $A$, whose cost is given by $2\times \text{nnz}(A)$.

The {compute time} of a sparse direct triangular solver on a many-core GPU
architecture such as from the NVIDIA cuSparse library is 
{much higher} than the SpMV.
For the NVIDIA V100 GPU architecture, the SpMV can now
achieve on the order of 50 -- 100 GigaFlops/sec in double precision
floating point arithmetic. 
When the number of GMRES$+$AMG iterations to achieve the same 
NRBE remains less than two times larger, then the case for
employing the mixed scheme on GPUs becomes rather compelling.
\begin{table}[htb!]
\centering
\begin{tabular}{|l|c|c|}
\hline
\textbf{Level}&$n_{\ell}$& $\text{nnz}(A)$ \\\hline
$\ell=$ 1 & 11498575 & 306389891  \\\hline
$\ell=$ 2 & 171074   & 6919886    \\\hline
$\ell=$ 3 & 19658    & 1256596     \\\hline
$\ell=$ 4 & 2018     & 126882      \\\hline
$\ell=$ 5 & 232      & 10070      \\\hline
$\ell=$ 6 & 29       & 605      \\\hline
$\ell=$ 7 & 5        & 25      \\\hline
\end{tabular}
\caption{Size and number of non-zeros for $A$, at each level, where
$n_\ell$ denotes the matrix size at level $\ell$. hypre-BoomerAMG.  
Aggressive coarsening. Matrix size $N = 11.49$M.}
\label{table:N14ka}
\end{table}

{
The compute time for the AMG $V$-cycle with GMRES SpMV is largely
determined by the memory bandwidth (BW) between the GPU main memory 
and the thread processors. For the AMD MI250X GPU employed on
Crusher and Frontier, the maximum achievable BW is 1.6 TeraBytes/sec.
Our performance estimates are based on a roof-line model for
the GPU and assumes we are in the bandwidth dominated regime.
We also assume that computation and data movement are not overlapped.
The transfer speed of sparse matrices (CSR format assumed)
between the GPU memory and the thread processors is the dominant cost.
More specifically, the $P/R$ prolongation-restriction matrices,
the coarse $A_c$, $L$ and $U$ factors, and the diagonal $D$ are read
at each level $i$ of the $V$-cycle hierarchy with the row sizes
and $nnz$ given in Table \ref{table:N14ka}. The transfer times
are estimated in the model given below, where it is assumed that
CSR matrix reads achieve $2/3$ of peak bandwidth.}
\begin{align*}
\text{Ac\_transfer} &= \text{nnz(i)} \times 8 \\
\text{L\_U\_transfer\_Jacobi} &= \text{nnz(i)} \times 8 \\
\text{P\_R\_transfer} &= 2 \times \text{nnz(i)} \times 8 \\
\text{D\_transfer\_Jacobi} &= \text{rows(i)} \times 8 \\
\text{matrix\_transfer\_cost} &= \text{nnz(1)} \times 8 \\
\text{A\_transfer} &= \text{matrix\_transfer\_cost} \\
\text{ILU}_{\text{Jacobi}}(\text{i}) &= \text{Ac\_transfer} + \text{P\_R\_transfer} + \text{D\_transfer\_Jacobi} 
\end{align*}
Assuming $N = 11 \times 10^6$ rows and $26$ nnz per row, and a computation efficiency of 200 GigaFlops/sec:
\begin{align*}
\text{FLOPs time} &\approx \frac{6 \times 40 \times 26 \times N}{200 \times 10^{9}} = 0.35 \text{ sec} \\
\text{Memory transfer time} &= \frac{2 \times \text{sum(ILU}_{\text{Jacobi}}) + \text{A\_transfer}}{1.6 \times 10^{12} \times 8 \times 40} = 3.05 \text{ sec} \\
\text{Total time} &= \text{FLOPs time} + \text{Memory transfer time} = 3.5\text{ sec} \\
\text{Measured} &: 3.7 \text{ sec}
\end{align*}

{
In practice, the above model is remarkably accurate. For the matrix
size $N = 11.49$M and $40$
GMRES+AMG iterations on a single GPU, we estimate a run time of 
$3.5$ sec, and measured, $3.7$ sec. In addition, the estimated
ratio
$\text{ILU}_{\text{GaussSeidel}} / \text{ILU}_{\text{Jacobi}} = 1.3$
is measured as $1.3$. Therefore the $V$-cycle with $\ell_1$-Jacobi is
30\% more efficient (faster) than when using Gauss-Seidel on the coarser
levels to achieve the same relative residual error tolerance of 1e$-6$.}

The compute times of the GMRES$+$AMG solver
in hypre with an incomplete $LDU$ smoother, and either direct or iterative 
solvers in the ILU smoother, are compared below for the NVIDIA V100 GPU.
The compute times for a single pressure solve
are given in Table \ref{tab:ilu-smoothers} for the $N=14186$ dimension
matrix. The ILU$(0)$ and ILUT smoothers are included for comparison. 
Both the CPU and GPU times are reported. In all cases, two 
Gauss-Seidel and one ILU sweep are employed. 
The solver time reported again corresponds
to when the relative residual decreases below 1$e-$5.

\begin{table}[htb]
\centering
\begin{tabular}{|l|c|c|c|c|c|} \hline
 &  Gauss-Seidel & Poly G-S & ILUT direct & ILUT iter & ILU$(0)$ iter \\ \hline
iterations &    7 &     9  &    7 &      7 &     5 \\ \hline
CPU (sec) &  0.037 & 0.025  & 0.025 & 0.035  & 0.038 \\ \hline
GPU (sec) &  0.021 & 0.0067  & 0.032 & 0.0065 & 0.0048 \\ \hline
\end{tabular}
\caption{ \label{tab:ilu-smoothers} GMRES$+$AMG compute time
on V100.
Gauss-Seidel, polynomial Gauss-Seidel, and ILU smoothers. 
PeleLM matrix dimension $N=14186$}
\end{table}

Consider the CPU compute times for a single solve.
The results indicate that the GMRES$+$AMG solver time using the
mixed $V$-cycle with an ILU$(0)$ smoother on the first level, versus a direct 
solver for the $L$ and $U$ systems, costs less than Gauss-Seidel
smoothing on all levels. The PMIS algorithm is employed
along with aggressive coarsening on the first $V$-cycle level.
One sweep of the Gauss-Seidel smoother is
employed in both configurations. The longer time is primarily
due to the higher number of Krylov iterations required to converge.
The ILUT smoother, on the finest level, with iterative solvers is the more
efficient approach on the GPU. Despite only ten (10) SpMV
products to solve the $L$ and 15 to solve the $U$ systems, 
the computational speed of the GPU for the
SpMV kernel is more than sufficient to overcome the cost of a direct solve. 
The observed solver time was 0.0065 with a performance model estimate of $0.0055$ sec
when $N=14186$ and requires twenty iterations 
to converge where aggressive coarsening with hypre reduces
the $\text{nnz}(A)$ per level \cite{Yang10}.

\begin{figure}[h!]
\centering
\includegraphics[height=0.35\textheight,width=0.65\textwidth]{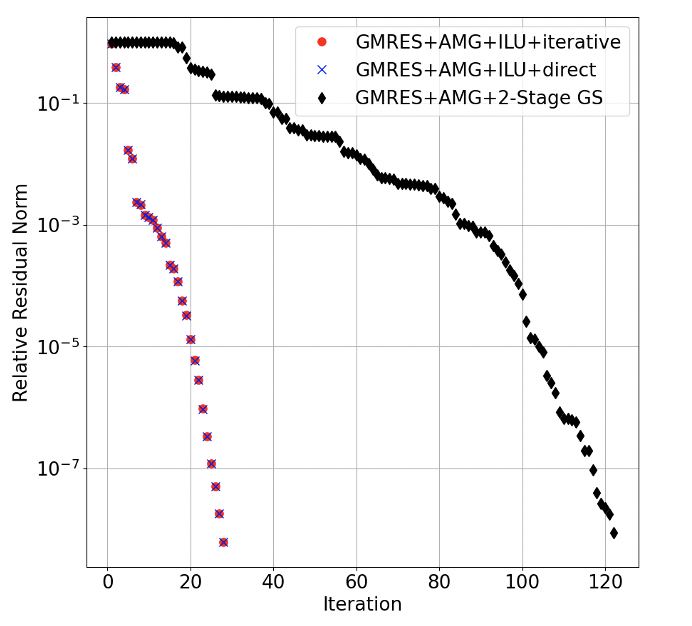}
\caption{\label{fig:11million-cvg} 
Convergence histories of (F)GMRES+AMG with polynomial  
Gauss-Seidel and ILU$(0)$ smoothers
using direct and iterative solves. $N = 11$M}
\end{figure}

\begin{figure}
\centering
\includegraphics[height=0.35\textheight,width=0.55\textwidth,]{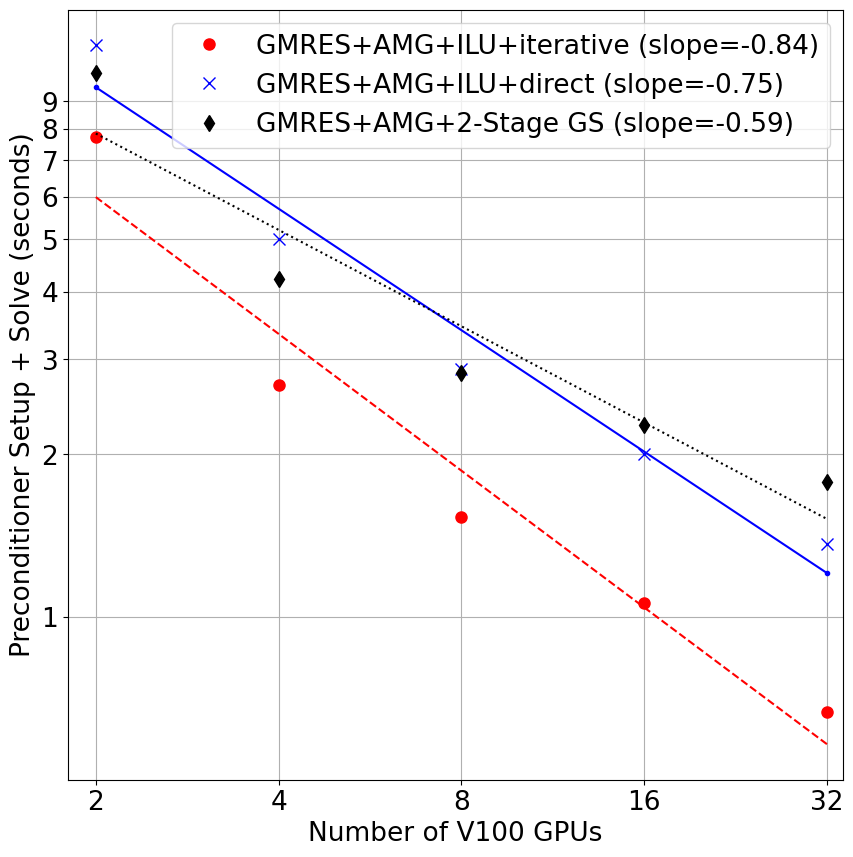}
\caption{\label{fig:11million} Strong-scaling of (F)GMRES+AMG with ILU$(0)$ smoothers
using direct and iterative solves. hypre-BoomerAMG GPU results on NREL Eagle.
Matrix size $N=11$M.}
\end{figure}

The compute times for a larger dimension problem where $N=1.4$M are
reported in Table \ref{tab:ilu-smoothers-2}. Here it was observed that the
ILU$(0)$ compute time on the GPU is lower than with the
polynomial Gauss-Seidel smoother.  
Most notably, the GPU compute time for ILU$(0)$ solves with fifteen (15) SpMV
for $U$ are two times faster than ILUT with direct  solves.
To further explore the parallel strong-scaling behaviour of the iterative
and direct solvers within the ILU smoothers, the GMRES+AMG solver was 
employed to solve a PeleLM linear system of dimensions $N=11$M.
The $LDU$ form of the factorization with row scaling was
again employed and twenty (20) SpMV provide sufficient smoothing
for this much larger problem. The linear system solver was tested 
on the NREL Eagle Supercomputer
configured with two NVIDIA Volta V100 GPUs per node.
Most notably, the solver with iterative Neumann scheme
achieves a faster solve time compared to the direct
solver as displayed in Figure \ref{fig:11million}.
The convergence histories of the GMRES+AMG solver with a 
polynomial  Gauss-Seidel and the ILU direct
and iterative smoothers are plotted in Figure
\ref{fig:11million-cvg}.

\begin{figure}[h!]
\centering
\includegraphics[height=0.35\textheight,width=0.65\textwidth]{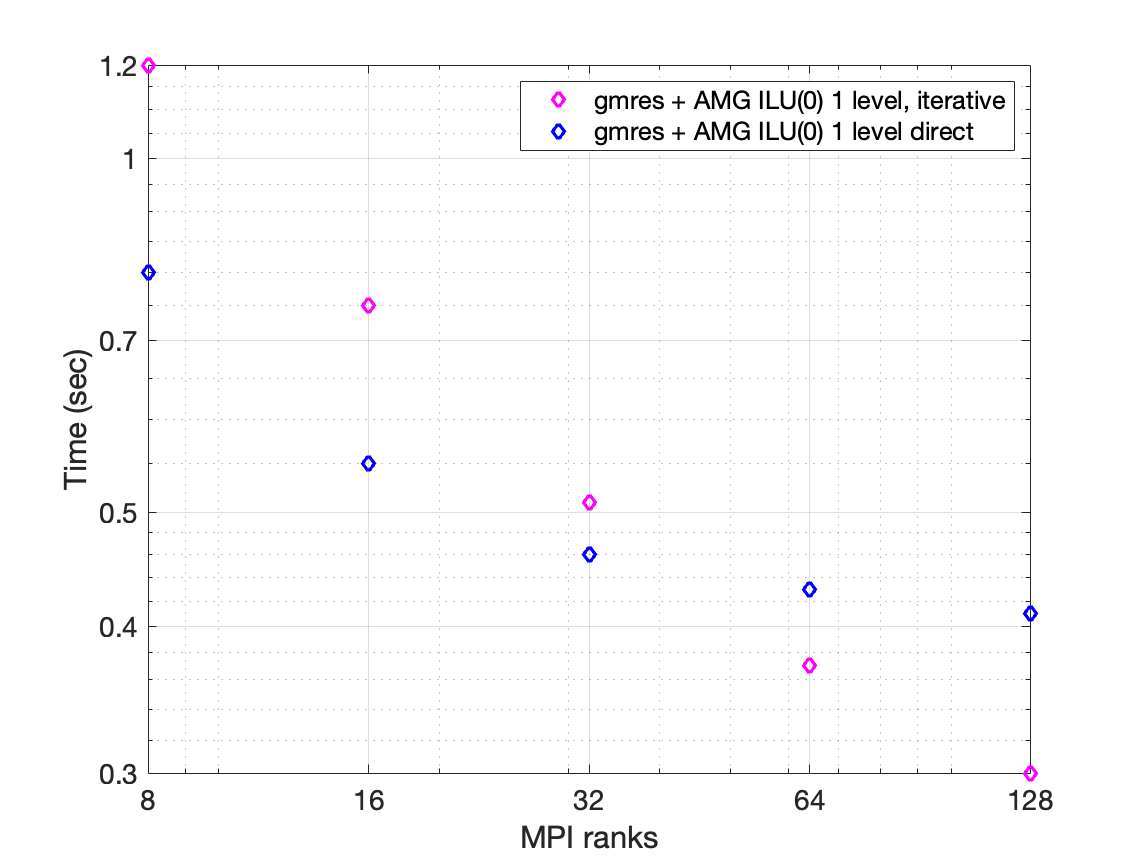}
\caption{\label{fig:11Mcrusher} 
hypre-BoomerAMG GPU results on ORNL Crusher with AND MI250X. PeleLM model. $N = 11$M}
\end{figure}

The larger problem was run on the ORNL Crusher supercomputer and the
compute times on up to 128 MPI ranks are displayed in
Figure \ref{fig:11Mcrusher}. The combined solve and set-up
times are plotted for direct and iterative triangular solvers.
For a low number of MPI ranks, the direct triangular solver
from AMD is more efficient than the iterative solves with SpMV.
However, as the number of ranks increases, the iterative 
triangular solver for the smoother leads to lower compute times.
The cross-over point for iterative solves occurs after thirty-two
GPUs where the compute time continues to drop significantly faster.
By comparing the times in Figures \ref{fig:11million} and 
\ref{fig:11Mcrusher}, it is observed that the AMD MI250X GPU
is roughly four times faster than the NVIDIA V100 for
this problem.

\begin{figure}
\centering
\includegraphics[height=0.35\textheight,width=0.65\textwidth]{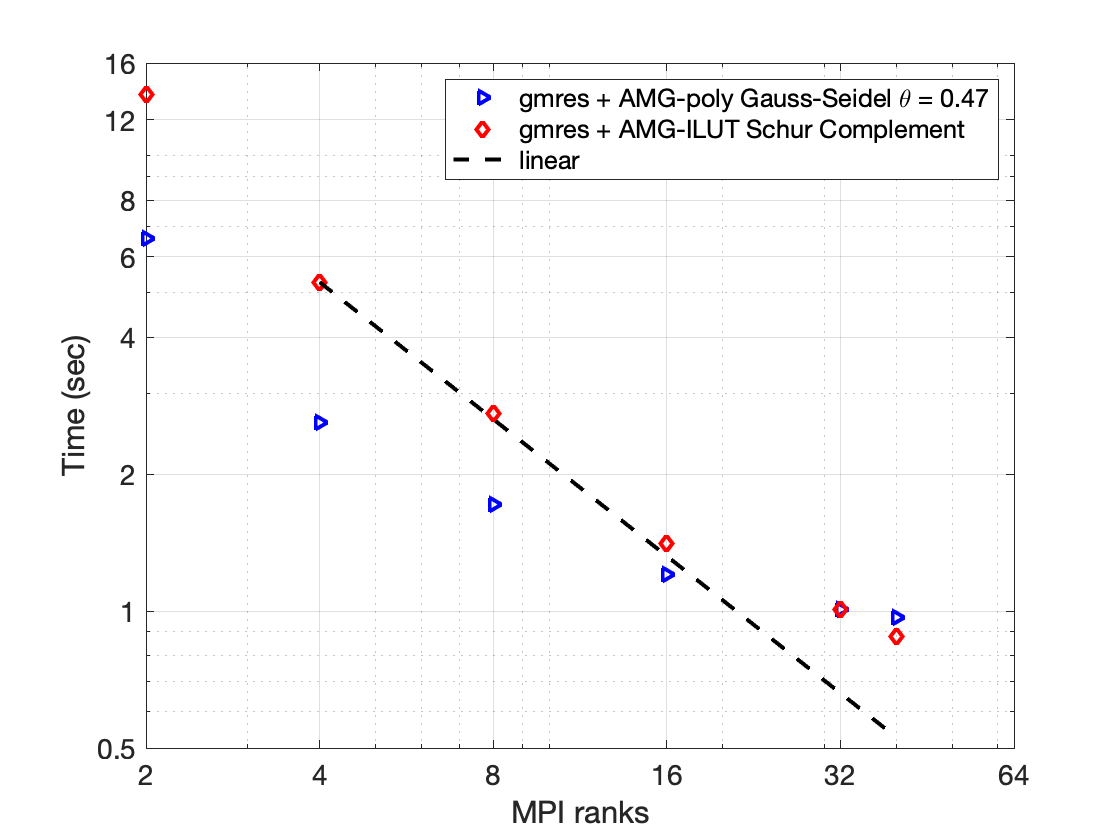}
\caption{\label{fig:23million-cvg}Nalu-Wind
NREL 5MW wind turbine mesh.
Strong-scaling of (F)GMRES+AMG with ILUT Schur Complement
using iterative solves versus polynomial Gauss-Seidel 
smoother. Matrix size $N=23$M}
\end{figure}

A strong-scaling study for the low resolution NREL 5 MegaWatt
single-turbine mesh is displayed in Figure \ref{fig:23million-cvg}. 
The matrix dimension for this problem is $N=23$M. The total
setup plus solve time is displayed for (F)GMRES+AMG executing
on the NREL Eagle supercomputer using two NVIDIA V100 GPUs
per node on up to 20 nodes or 40 GPUs. The solve time is
plotted for the polynomial Gauss-Seidel and ILUT Schur
complement smoothers. In the latter case, a single 
iteration of the iterative GMRES solver,
without residual computations, results
in a steeper decrease in the execution time
and improved performance. 
In addition, the number of GMRES+AMG solver iterations
to reach a relative residual tolerance of 1e$-5$
remains constant at eleven (11) as the number of
compute nodes is increased.

\begin{table}[htb]
\centering
\begin{tabular}{|l|c|c|c|c|c|} \hline
 &  Gauss-Seidel & Poly G-S & ILUT direct & ILUT iter & ILU$(0)$ iter
 \\ \hline
 iterations &    7 &       9  &     8 &     8 &      4 \\ \hline
CPU (sec) &    9.2 &     9.6  &   4.3 &   6.9 &    6.8 \\ \hline
GPU (sec) &   0.29 &   0.055  & 0.098 & 0.058 &  0.042 \\ \hline
\end{tabular}
\caption{ \label{tab:ilu-smoothers-2} GMRES$+$AMG compute time
on V100. Gauss-Seidel, polynomial Gauss-Seidel, and ILU smoothers. 
PeleLM matrix dimension $N=1.4$M.}
\end{table}

\section{Conclusions} \label{sec:concl}

For the highly ill-conditioned PeleLM nodal pressure projection linear systems,
the standard Jacobi and Gauss-Seidel smoothers are less effective for
reducing the error at each level of the AMG
$V$-cycle and may result in very large iteration counts or fail to 
converge for the preconditioned  Krylov solver. 
Jomo et al.~\cite{Jomo2021} compare PCG+AMG
with Jacobi and Gauss-Seidel smoothers with GMRES preconditioned
by ILU. However, they did not investigate ILU as a smoother for AMG.
We proposed a novel approach for the solution
of sparse triangular systems for the
$L$ and $U$ factors of an ILU smoother for AMG.
Previous work by H. Anzt, and E. Chow demonstrated that these
factors are highly non-normal,
even after appropriate re-ordering and scaling of 
the linear system $Ax = b$ and,
Jacobi iterations may diverge.
In order to mitigate the effects of a high degree of non-normality, 
as measured by Henrici's metric, either a row or
row/column scaling is applied to the $U$ factor during the
set-up phase of the AMG $V$-cycle. A finite Neumann series
multiplied by a vector is then computed, which is equivalent to
a Richardson iteration. Our results demonstrate that
several orders of magnitude reduction in the departure
from normality ${\rm dep}(U)$ is possible, thus leading to 
robust convergence of GMRES+AMG.

In order to further improve the efficiency of
the PeleLM (F)GMRES$+$AMG pressure solver on 
many-core GPU architectures, 
AMG $V$-cycles with ILU smoothing
on the finest levels are combined with a polynomial smoother
on the remaining coarse levels. It was found that the convergence rates
for mixed AMG are almost identical to using ILU on all
levels, thus leading to significant cost reductions.
For a large problem $N=11$M
solved on the NREL Eagle supercomputer the iterative 
solve for $LDU$ with row scaling led to a 5X
speed-up over the direct solve within the GMRES+AMG 
$V$-cycles. Furthermore, the strong scaling curve for the
solver run time is close to linear.

For the Nalu-Wind pressure continuity equation, 
an ILUT Schur complement smoother
with iterative solves on the local block diagonal
systems was applied. Pressure linear systems from NREL 5 MegaWatt 
reference turbine simulations were employed to assess numerical
accuracy and performance. The linear solver exhibits improved
parallel strong-scaling characteristics with this new
smoother and maintains a constant GMRES iteration count
when the number of GPU compute nodes increases. 
By omitting the residual computations from
the single iteration of the GMRES-Schur solver, 
to overall execution time is reduced.
Our future plans include implementing the fixed-point iteration 
algorithms of Chow 
\cite{ChowPatel2015} to compute the ILU
factorization on GPUs. The Schwarz preconditioners
described in Prenter \cite{Prenter2020} and Jomo et al.
\cite{Jomo2021} could also be adapted to hypre for PeleLM.
The solver has also been incorporated into
the MFIX-Exa CFD-DEM model \cite{Musser2022}
for carbon capture and chemically reacting fluid flows.
The solver leads to at least 5X improvement in
the computational speed on GPUs.

\section*{Acknowledgments}
This work was authored in part by the National Renewable Energy Laboratory, 
operated by 
Alliance for Sustainable Energy, LLC, 
for the U.S. Department of Energy (DOE) under Contract 
No. DE-AC36-08GO28308. Funding was provided by the Exascale Computing Project 
(17-SC-20-SC), a collaborative
effort of two U.S. Department of Energy organizations (ASCR and 
the NNSA). 
{The data that support the findings of this study are available 
from the corresponding author upon reasonable request.
We thank the anonymous reviewers for their suggestions, 
which led to several valuable improvements.}
%

\bibliography{paper-ILUTPforAMG}%


\end{document}